\documentclass[12points]{amsart}

\usepackage{amsmath,amssymb,amsfonts,graphics,epsfig}
\usepackage[matrix,arrow,curve]{xy}

\newcommand{\DP}[2]{{\frac{\partial #1}{\partial #2}}}

\newcommand{\tx}[1]{\text{\rm #1}}
\newcommand{\CC}{\mathbb{C}}

\newcommand{\RR}{\mathbb{R}}   
\newcommand{\Ss}{\mathbb{S}} 
 
\newcommand{\ZZ}{\mathbb{Z}}      

\newcommand{\demo}{\noindent {\it \small Proof. }} 
 
\newtheorem{defi}{Definition}

\newtheorem{prop}{Proposition} 
\newtheorem{lem}{Lemma}

\def\mc{\mathcal}

\title[The K\"unneth formula in Floer homology]
{The K\"unneth formula in Floer homology for manifolds with restricted 
contact type boundary}
   
\author{{Alexandru Oancea}}

\date{24 August 2005} 

\begin{document}

\maketitle

\vspace{-.8cm}

\begin{center} 
Department of Mathematics, ETHZ,

R\"amistrasse 101, 8092 Z\"urich (CH).

Email: {\tt
oancea\,@\,math.ethz.ch}

\end{center}

\medskip

\begin{abstract} 

We prove the K\"unneth formula in Floer
(co)homology for manifolds with restricted contact type boundary. We
use Viterbo's definition of Floer homology, involving the
symplectic completion by adding a positive cone over the boundary. 
The K\"unneth formula implies 
the vanishing of Floer (co)homology 
for subcritical Stein manifolds. Other applications include the
Weinstein conjecture 
in certain product manifolds, obstructions to exact 
Lagrangian embeddings, existence of holomorphic curves with Lagrangian
boundary condition, as well as symplectic capacities. 

\end{abstract}

\renewcommand{\thefootnote}{\fnsymbol{footnote}}
\setcounter{footnote}{0}
\footnotetext{ 

The author is currently supported by the 
Forschungsinstitut f\"ur Mathematik at ETH Z\"urich.
}

\renewcommand{\thefootnote}{\arabic{footnote}}
\setcounter{footnote}{0}

\section{Introduction}

The present paper is concerned with the Floer homology groups 
$FH_*(M)$ of a compact symplectic manifold $(M, \, \omega)$ 
with contact type
boundary, as well as with their cohomological dual analogues
$FH^*(M)$. The latter were defined by Viterbo in \cite{functors1}
and are invariants that take into account the topology of the underlying
manifold {\it and}, through an algebraic limit process, all closed
characteristics on $\partial M$. Their definition is closely related
to the Symplectic homology groups of Floer, Hofer, Cieliebak and
Wysocki \cite{FH,CFH,FHW,CFHW,Cieliebak handles}.
%, as well as to other recent
%constructions of Biran, Polterovich and Salamon~\cite{BPS} or
%Frauenfelder and Schlenk~\cite{FrSc}. 

Throughout this paper we will assume that $\omega$ is exact, 
and in particular\break 
$\langle\omega,\,
\pi_2(M)\rangle=0$. This last condition will be referred to as {\it symplectic 
  asphericity}. 
%This condition holds for exact manifolds and 
%in particular for manifolds with restricted contact
%type boundary.  
The 
groups $FH_*(M)$ are invariant with respect to deformations of the
symplectic form $\omega$ that preserve the contact type character of
the boundary and the condition $\langle \omega,\, \pi_2(M)
\rangle=0$. The groups $FH_*(M)$ actually depend 
only on the {\it symplectic
  completion} $\widehat M$ of $M$. The manifold $\widehat M$ is
obtained by gluing a positive cone along the boundary $\partial M$
and carries a symplectic form $\widehat \omega$ which is canonically
determined by $\omega$ and the conformal vector field on $M$. We shall
often write $FH_*(\widehat M)$ instead of $FH_*(M)$. 
The grading on $FH_*(\widehat M)$ is given by minus the Conley-Zehnder
index modulo $2\nu $, with $\nu$ the minimal Chern number of $M$. There
exist canonical maps 
$$H_{n+*}(M,\, \partial M) \stackrel {c_*} \longrightarrow FH_*(\widehat M),$$
$$FH^*(\widehat M) \stackrel {c^*} \longrightarrow H^{n+*}(M, \, \partial M)$$
which shift the grading by $n=\frac 1 2 \dim M$. 
%The failure 
%of these morphisms to be injective, respectively 
%surjective has deep implications~\cite{functors1} on the
%symplectic geometry of $M$. Our main result is the following. 

%A simplified version of our main result is the following (see Theorem 
%\ref{theoreme Kunneth} 
%%and \ref{theoreme Kunneth cohomologique} 
%for the full statement). 

\medskip 

\noindent {\bf Theorem A (K\"unneth formula).} 
{\it Let $(M^{2m}, \, \omega)$ and $(N^{2n}, \,
  \sigma)$ be compact symplectic
manifolds with {\rm restricted} contact type boundary. 
Denote the minimal Chern numbers
 of $M$, $N$ and $M\times N$ by $\nu_M$, $\nu_N$ and $\nu_{M\times N}
 = \tx{gcd}\,(\nu_M,\, \nu_N)$ respectively.    
\renewcommand{\theenumi}{\alph{enumi}} 
\begin{enumerate}    
 \item For any ring $A$ of coefficients there exists a short exact sequence  
which splits noncanonically 
\begin{equation} \label{suite exacte Kunneth Floer}  
{\scriptsize
\xymatrix  
@C=30pt  
@R=20pt@W=1pt@H=1pt  
{
 {\bigoplus_{\widehat r+\widehat s = k} }
FH_{\widehat r}(M, \, \omega) \otimes    
  FH_{\widehat s}(N, \, \sigma) \ \ \ar@{>->}[r] & 
FH_k(M \times N, \, \omega \oplus \sigma )    
  \ar@{->>}[d]    \\    
  & {\bigoplus_{\widehat r+\widehat s = k-1}} 
  \tx{Tor}_1^A \big(
  FH_{\widehat r}(M, \, \omega), \    
  FH_{\widehat s}(N, \, \sigma) \big) 
}     
}
\end{equation}

The morphism $c_*$ induces a morphism of exact sequences 
whose source is the K\"unneth
exact sequence of the product $(M,\, \partial M) \times (N,\,
\partial N)$ and whose target is~(\ref{suite exacte Kunneth Floer}).

\item For any field $\mathbb{K}$ of coefficients there is an isomorphism   
\begin{equation} \label{Kunneth egal Floer coh corps}  
{\scriptsize
\xymatrix  
@C=30pt  
@R=20pt@W=1pt@H=1pt  
{ \bigoplus_{\widehat r+\widehat s = k} FH^{\widehat r}(M, \, \omega) 
\otimes    
  FH^{\widehat s}(N, \, \sigma) \ar[r]^{\qquad \sim}  
& FH^k(M \times N, \, \omega \oplus 
  \sigma ) \ ,    
}     
}
\end{equation} 

The morphism $c^*$ induces a commutative diagram with respect to the
K\"unneth isomorphism in cohomology for $(M,\, \partial
M)\times (N,\, \partial N)$.
\end{enumerate} 

%In the above notation we have $k\in \ZZ/2\nu_{M\times N}\ZZ$, $0 \le
%r \le 2\nu _M -1$, 
%$0\le s \le 2\nu_N -1$ and the $\widehat{\ }$ symbol associates
%to an integer its class in the corresponding $\ZZ/2\nu\ZZ$ ring.

%Then Floer cohomology with coefficients in a field satisfies 
%\begin{equation} \label{iso K coh}
%FH^*(\widehat M) \otimes
%FH^*(\widehat N) \simeq FH^*(\widehat M\times \widehat N) \ .
%\end{equation}
%The product $\widehat M\times \widehat N$ 
%is endowed with the product symplectic form
%$\widehat \omega \oplus \widehat \sigma$. The above
%isomorphism is compatible with the maps $c^*$ and the usual K\"unneth
%isomorphism 
%\begin{equation*}
%  H^{*}(M, \, \partial M) \otimes H^{*}(N, \, \partial N)\simeq H^*(M
%  \times N, \, \partial (M \times  N) \ .
%\end{equation*}

%Moreover, there is a commutative diagram 

%\begin{equation} \label{c d K coh}
%\xymatrix  
%@C=13pt  
%@R=30pt@W=1pt@H=1pt  
%{  FH^*(\widehat M) \otimes    
%  FH^*(\widehat N) \ar[rr]^{\quad \quad \ \ \sim} \ar[d]_{c^*\otimes c^*}   
%& & FH^*(\widehat M \times \widehat N) \ar[d]^{c^*}  \\ 
% H^{*}(M, \, \partial M) \otimes    
%  H^{*}(N, \, \partial N) \ar[rr]^{\quad  \sim} & &  
%H^*(M \times N, \, \partial (M \times  
%N) )  
%} 
%\end{equation}
}
 
\medskip 

In the above notation we have $k\in \ZZ/2\nu_{M\times N}\ZZ$, $0 \le
r \le 2\nu _M -1$, 
$0\le s \le 2\nu_N -1$ and the $\widehat{\ }$ symbol associates
to an integer its class in the corresponding $\ZZ/2\nu\ZZ$ ring. The
reader can consult~\cite[VI.12.16]{D} for a construction of the 
K\"unneth exact sequence in singular homology.

\medskip

The algebraic properties of the map $c^*$ strongly influence
the symplectic topology of the underlying manifold. Our applications
are based on the following theorem, which summarizes part of
the results in \cite{functors1}. 

\medskip 

\noindent {\bf Theorem (Viterbo \cite{functors1}).} {\it Let $(M^{2m},\, \omega)$ be a
  manifold with contact type boundary such that $\langle \omega,\,
  \pi_2(M)\rangle=0$. Assume the map $c^*:FH^*(M)
  \longrightarrow H^{2m}(M, \, \partial M)$ is {\rm not} surjective. 
  Then the following hold.
  \renewcommand{\theenumi}{\alph{enumi}}
  \begin{enumerate}
   \item The same is true for any hypersurface of restricted contact
   type $\Sigma \subset M$ bounding a compact region;
   \item Any hypersurface of contact type $\Sigma \subset M$ bounding
   a compact region carries a closed characteristic (Weinstein conjecture);
   \item There is no exact Lagrangian embedding $L \subset M$ (here
   $M$ is assumed to be exact by definition);
   \item For any Lagrangian embedding $L \subset M$ there is a loop on
   $L$ which is contractible in $M$, has strictly positive area and
   whose Maslov number is at most equal to $m+1$;
   \item For any Lagrangian embedding $L \subset M$ and any compatible
   almost complex structure $J$ there is a nonconstant 
   $J$-holomorphic curve $S$
   (of unknown genus)
   with non-empty boundary $\partial S \subset L$.   
  \end{enumerate}

}

\medskip 

Viterbo
\cite{functors1} introduces
the following definition, whose interest is obvious in the light of the
above theorem.

\begin{defi} \tx{(Viterbo)}
    A symplectic manifold $(M^{2m}, \, \omega)$ which verifies\break $\langle
  \omega, \, \pi_2(M) \rangle = 0$ is said to satisfy the Strong
  Algebraic Weinstein Conjecture (SAWC) if the composed morphism below is not surjective
  $$FH^*(M) \stackrel{c^*}\longrightarrow H^*(M, \, \partial M)
  \stackrel{\tx{pr}}\longrightarrow 
  H^{2m}(M, \, \partial M) \ .$$ 
\end{defi}

We shall still denote the composed morphism by $c^*$. 
%The condition $\langle [\omega], \, \pi_2(M) \rangle = 0$ 
%greatly simplifies
%the analysis involved in the construction of Floer homology as there can be no
%bubbling-off of pseudo-holomorphic spheres in $M$, these having
%positive area. 
Theorem A now implies that the
property of satisfying the SAWC is stable under products, 
with all the geometric consequences listed above. 

\medskip 
 
\noindent {\bf Theorem B.} {\it Let $(M, \, \omega)$ 
  be a symplectic manifold
  with restricted contact type boundary satisfying the SAWC. 
  Let  $(N, \, \sigma)$ be an arbitrary symplectic manifold with
  restricted contact type boundary. The product $(\widehat
  M \times \widehat N, \, \widehat \omega \oplus \widehat \sigma)$
  satisfies the SAWC and assertions (a) to (e) in the above theorem of
  Viterbo hold. In particular, the Weinstein conjecture holds
  and there is no exact Lagrangian embedding in $\widehat M\times \widehat N$. 
}

\medskip 

The previous result can be 
 applied for subcritical
 Stein manifolds of finite type. 
These are complex manifolds $\widehat
 M$ which
 admit proper and bounded from below plurisubharmonic Morse functions 
 with only a finite number of critical points, all of index strictly
 less than $\frac 1 2 \dim _\RR \widehat M$~\cite{Eli-psh}. 
They
 satisfy the SAWC as proved by 
 Viterbo~\cite{functors1}. 
%Note however that the property of
% having a restricted contact type boundary is strictly more general than that
% of being a Stein domain. 
%Taking the product with a  
% Stein manifold can
% therefore produce restricted contact type manifolds which are not
% Stein anymore. 
Cieliebak~\cite{Cieliebak handles} has proved that their Floer homology
actually vanishes. We can recover this
through Theorem A by using another of his results~\cite{C}, 
namely that every such manifold is 
Stein deformation equivalent to a split one $(V \times \CC, \, \omega
\oplus \omega_{\tx{std}})$. This can
be seen as an extension of the classical vanishing result
$FH_*(\CC^\ell)=0$, $\ell \ge 1$~\cite{FHW}. 

\medskip 

\noindent {\bf Theorem C (Cieliebak~\cite{Cieliebak handles}).} {\it
  Let $\widehat M$ be a  
subcritical Stein manifold of finite type. Its Floer homology vanishes} 
$$FH_*(\widehat M) = 0 \ .$$

%To the best of my knowledge there are only three other explicit
%computations of groups $FH^*(\widehat N)$, 
%namely $FH^*(\CC^n)=0$, $n \ge 1$ \cite{FHW},  
%$FH^*(T^*N)
%\simeq H^*(\Lambda_0 N)$, with
%$N$ a closed manifold and $\Lambda_0 N$   
%the space of contractible loops \cite{cot bdls Viterbo,Salamon
%  Weber}, and $FH^*(\mathcal L)=0$, 
%with $\mathcal L$ a line bundle with negative Chern class over
%a symplectically aspherical manifold \cite{teza mea}.

%However, our applications will all concern manifolds for
%which $c^*$ is not surjective, and this has already been proved to be
%the case for
%subcritical Stein manifolds by Viterbo \cite{functors1}. 

\medskip

The paper is organized as follows. In Section~\ref{constructions} we
state the relevant definitions and explain the main properties of
the invariant $FH_*$. 
%We state the homological couterpart of Theorem A
%in Section~\ref{KuKu}, where we also give a sketch of the proof. 
Section~\ref{la preuve
  de Kunneth} contains the proof of Theorem A. 
The proofs of Theorems B and C, together with other applications, 
are gathered in Section~\ref{appli}.

\medskip 
 
Let us point out where the
difficulty lies in the proof of Theorem A. 
Floer homology is 
defined on closed manifolds for any Hamiltonian satisfying some
generic nondegeneracy 
condition, and this condition is {\it stable} under sums $H(x)+K(y)$ on
products $M \times N \ni (x, \, y)$. This trivially implies (with
field coefficients) a
K\"unneth formula of the type $FH_*(M\times N; \, H+K, \, J_1\oplus
J_2) \simeq FH_*(M; \, H, \, J_1) \otimes FH_*(N; \, K, \, J_2)$. On
the other hand, Floer homology for manifolds with contact type boundary
is defined using Hamiltonians with a rigid behaviour at infinity and
involves an algebraic limit construction.  
This class of
Hamiltonians is {\it not stable} under the sum operation $H(x)+K(y)$
on $M \times N$. One may still define Floer homology groups
$FH_*(M\times N; \, H+K, \, J_1\oplus J_2)$,  but
%, as $H+K$ no more
%belongs to the admissible class, 
the resulting homology might well be
different, in the limit, 
from $FH_*(\widehat M\times \widehat N)$. The whole point of
the proof is to show that this is not the case.

\medskip 

This paper is the first of a series studying the Floer homology of 
symplectic fibrations with contact type boundary. It treats
trivial fibrations with open fiber and base. A spectral
sequence of Leray-Serre type for symplectic
fibrations with closed base and open fiber is constructed in~\cite{LS}.

\medskip 

{\it Acknowledgements.} This work is part of my Ph.D. 
thesis, which I completed under the guidance of Claude Viterbo. Without his
inspired 
support this could not have come to being. I am grateful to 
Yasha Eliashberg,  
Dietmar Salamon, Paul Seidel, Jean-Claude Sikorav and Ivan Smith 
for their help and suggestions. I also 
thank the referee for pointing out 
errors in the
initial proofs of Theorem C and Proposition~\ref{prop:symplectic
  capacities}.
%, as well as the reference~\cite{Cieliebak handles} 

During the various stages of preparation of this work I was supported 
by the following institutions: Laboratoire de
Math\'ematiques, Universit\'e Paris Sud~; 
Centre de Math\'ematiques de l'\'Ecole Polytechnique~;
\'Ecole Normale Sup\'erieure de Lyon~; 
Departement Mathematik, ETH.

\section{Definition of Floer homology} \label{constructions}

 Floer homology has been first defined by A. Floer for closed
 manifolds in a  series of papers \cite{F1,F2} which
 proved Arnold's conjecture for a large class of
 symplectic manifolds including {\it symplectically aspherical} ones. 
%These
% are symplectic manifolds $(M, \, \omega)$ such that $\langle
% [\omega], \, \pi_2(M)\rangle = 0$.  
%
%\footnote{The
% currently used  terminology calls ``symplectically aspherical'' a
% manifold such that $\langle [\omega], \, \pi_2(M) \rangle =0$. I have
% decided to 
% use this name in the stronger sense defined in the text in order to avoid
% cumbersome formulations like ``strongly symplectically aspherical''.}  
% The first Chern class is understood with respect to an
% $\omega$-compatible almost complex structure $J$ on $TM$, where 
% {\it $\omega$-compatible} means  that $\omega(\cdot, \, J\cdot)$ is
%  a Riemannian metric. It is a classical result \cite{G} 
%that the set of all
% such almost complex structures is non empty and contractible, hence
% the first Chern class is independent of the particular choice of $J$.  
In this situation Floer's construction can be summarized as follows. Consider a
periodic 
time-dependent Hamiltonian $H : \Ss^1 \times M \longrightarrow \RR$
with Hamiltonian vector field $X^t_H$ defined by $\iota _{X^t_H}\omega
= dH(t, \cdot)$. The
associated {\it action functional}  
$$A_H : C^\infty _{\tx{contr}} (M) \longrightarrow \RR \ ,$$
$$\gamma \longmapsto -\int_{D^2}\bar{\gamma}^* \omega -
\int_{\Ss^1}H(t, \, \gamma(t)) dt $$
is defined on the space of smooth contractible loops 
$$C^\infty _{\tx{contr}} (M) = \big\{ \gamma:\Ss^1 \longrightarrow M \
: \ \exists \ \tx{smooth } \bar{\gamma}: D^2 \longrightarrow M, \
\bar{\gamma}|_{\Ss^1} = \gamma \big\} \ .$$
The critical points of $A_H$ are precisely the $1$-periodic solutions of
$\dot{\gamma}=X_H^t(\gamma(t))$, and we denote the corresponding set
by $\mc P(H)$. We suppose that  the elements of $\mc P(H)$ 
are nondegenerate i.e. the time one return map has no
eigenvalue equal to $1$. Each such periodic orbit $\gamma$ has a
$\ZZ/2\nu \ZZ$-valued  
Conley-Zehnder index $i_{CZ}(\gamma)$ (see \cite{RS1} for
the definition), where $\nu $ is the {\it minimal Chern
number}. The latter is defined by $\langle c_1, \, \pi_2(M) \rangle =
\nu \ZZ$ and by the convention $\nu = \infty$ if $\langle c_1, \, \pi
_2(M) \rangle =0$. 
Let also $J$ be a compatible almost complex
structure. The homological Floer complex is defined as 

$$FC_k(H, \, J) = \bigoplus_{\scriptsize 
\begin{array}{c} \gamma \in \mc P (H) \\
    i_{CZ}(\gamma) = -k \mod \ 2\nu 
 \end{array} } \ZZ \langle \gamma \rangle \ , $$

$$\delta : FC_k(H, \, J) \longrightarrow FC_{k-1}(H, \, J) \ ,$$

\begin{equation} \label{the differential} 
\delta \langle \gamma \rangle = \sum_{\scriptsize 
\begin{array}{c} 
\gamma' \in \mc P (H) \\ i_{CZ}(\gamma') = -k + 1 \mod \ 2\nu  
\end{array}} \# \mc M
( \gamma, \, \gamma'; \, H , \, J)/\RR \ . 
\end{equation} 

Here $\mc M
( \gamma, \, \gamma'; \, H , \, J)$ denotes the space of trajectories
for the negative $L^2$ gradient of $A_H$ with respect to the metric
$\omega(\cdot, \, J\cdot)$, running from $\gamma$ to $\gamma'$: 
\begin{equation} \label{Fl eq}
\mc M
( \gamma, \, \gamma'; \, H , \, J) = \Bigg\{ 
u: \RR \times \Ss^1 \longrightarrow M \, : \, 
\begin{array}{c} u_s + J\circ u \cdot u_t - \nabla H(t, \, u) =0 \\
u(s, \cdot) \longrightarrow \gamma, \ s \longrightarrow -\infty \\ 
u(s, \cdot) \longrightarrow \gamma' , \ s\longrightarrow +\infty 
\end{array}
\Bigg\}.
\end{equation} 

The additive group $\RR$ acts on $\mc M(\gamma , \, \gamma'; H , \,
J)$ by translations in the $s$ variable, while the symbol $\# \mc
M(\gamma, \, \gamma'; \, H, \, J)/\RR$ stands for an algebraic count of
the elements of $\mc M(\gamma, \, \gamma'; \, H, \, J)/\RR$. We note
  that it is possible to choose any coefficient ring instead of $\ZZ$ 
  once the sign assignment procedure is available.

The crucial statements of the theory are listed below. The main point is
that the symplectic asphericity condition prevents the loss of
compactness for Floer trajectories by bubbling-off of nonconstant 
$J$-holomorphic
spheres: the latter simply cannot exist. 
\renewcommand{\theenumi}{\alph{enumi}}
\begin{enumerate}
\item under the nondegeneracy assumption on the $1$-periodic orbits
  of $H$ and for a generic choice of $J$, 
  the space $\mc M(\gamma , \, \gamma'; H , \, J)$ is a
  manifold of dimension $i_{CZ}(\gamma') - i_{CZ}(\gamma) = 
  -i_{CZ}(\gamma) - ( -i_{CZ}(\gamma'))$. If this difference 
 is equal to $1$, 
  the space $\mc M(\gamma, \, \gamma'; H , \, J)/\RR$ consists of a
  finite number of points. Moreover, there is a consistent choice of
  signs for these points with respect to which $\delta \circ \delta =
  0$; 
\item the Floer homology groups $FH_*(H, \, J)$ are independent of $H$
  and $J$. More precisely, for any two pairs $(H_0, \, J_0)$ and
  $(H_1, \, J_1)$ there is a generic choice of a smooth 
  homotopy $(H_s, \, J_s)$, $s\in
  \RR$ with $(H_s, \, J_s)\equiv (H_0, \, J_0)$, $s\le 0$, $(H_s,
  \, J_s) \equiv (H_1, \, J_1)$, $s\ge 1$ defining a
  map 
\begin{equation} \label{cont morph}
  \sigma ^{(H_0, \, J_0)} _{(H_1, \, J_1)} : FC_k (H_0, \, J_0)
  \longrightarrow FC_k(H_1, \, J_1) \ ,
\end{equation}
  $$\sigma \langle \gamma \rangle = \sum_{\scriptsize \begin{array}{c}
  \gamma' \in \mc P (H_1) \\ i_{CZ}(\gamma')=i_{CZ}(\gamma) = -k \mod
  \ 2\nu 
  \end{array}} \# \mc M (\gamma, \, \gamma' ; \, H_s, \, J_s) \ .$$
  The notation $\mc M (\gamma, \, \gamma'; \, H_s, \, J_s ) $ stands
  for the space of solutions of the equation 
\begin{equation} \label{par Fl eq}
u_s + J(s, \, u(s, \, t) ) u_t - \nabla H(s, \, t, \, u(s, \, t))
  =0 \ ,\end{equation} 
  which run from $\gamma$ to $\gamma'$. The map $\sigma ^{(H_0, \,
  J_0)} _{(H_1, \, J_1)} $ induces an isomorphism in cohomology
  and this isomorphism is independent of the choice of the
  homotopy. We shall call it in the sequel the ``continuation morphism''. 

\item if $H$ is time independent, Morse and sufficiently small in some
  $C^2$ norm, then the $1$-periodic orbits of $X_H$ are the critical
  points of $H$ and the Morse and Conley-Zehnder indices satisfy
  $i_{\tx{Morse}}(\gamma) = m + (-i_{CZ}^0(\gamma))$, $m= \frac 1 2 \dim
  M$. Here $i_{CZ}^0(\gamma)$ is the Conley-Zehnder computed with
  respect to the trivial filling disc. 
  Moreover, the Floer trajectories running between points with
  index difference equal to $1$ are independent of the $t$ variable
  and the Floer complex is equal, modulo a shift in the grading, with
  the Morse complex of $H$.
 
We infer that for any regular pair $(H, \,
  J)$ we have 
$$\displaystyle FH_k(H, \, J) \simeq \bigoplus_{l \ \equiv \ k \mod \ 2\nu }
  H_{l+m}(M) \ , \ k \in \ZZ/2\nu \ZZ \ . $$ 
\end{enumerate}

\medskip 

{\bf Remark.} An analogous construction yields a {\it
  cohomological complex} by considering in (\ref{the differential}) the 
space of
  trajectories $\mc M(\gamma', \, \gamma; \, H, \, J)$. 

\medskip 

We explain now how the above ideas 
can be adapted in order to construct a symplectic
invariant for manifolds with contact type boundary. We
follow~\cite{functors1}, but similar constructions can be found
in~\cite{FH,CFH,Cieliebak handles} (see~\cite{my survey} for a
survey).  
%For a comparison
%of the various constructions in \cite{FH,CFH,functors1} one can
%consult the survey \cite{my survey}. 

\begin{defi} 
 A symplectic manifold  $(M, \, \omega)$ is said to have a {\rm
 contact type boundary} if there is a vector field $X$ defined in a
 neighbourhood of $\partial M$, pointing outwards and transverse to
 $\partial M$, which  satisfies 
 $$L_X\omega = \omega \ .$$
We say that  $M$ has {\rm restricted contact type boundary} if 
$X$ is defined globally.
\end{defi}

 We call $X$ and $\lambda = \iota(X) \omega$ the {\it Liouville vector
   field} and the {\it  Liouville
 form} respectively, with $d\lambda=\omega$. 
We define the {\rm Reeb vector field}
$X_{\tx{Reeb}}$ as the
generator of $\ker \omega|_{T\partial M}$ normalized by
$\lambda(X_{\tx{Reeb}})= 1$. An integral curve of $X_{\tx{Reeb}}$ is
called a {\it characteristic}. We have $\varphi_t^*
\omega = e^t\omega$ and this implies that 
a neighbourhood of $\partial M$ is foliated
by the hypersurfaces $\big( \varphi_t(\partial M) \big) _{-\epsilon
  \le t \le 0}$, whose characteristic dynamics are conjugate.

\medskip 

The Floer homology groups $FH_*(M)$ of a  manifold with contact type 
boundary are an invariant that takes
into account:
\begin{itemize} 
\item the dynamics  
on the boundary by ``counting'' characteristics of arbitrary period;
\item  the interior topology of
$M$ by ``counting'' interior $1$-periodic orbits of Hamiltonians.
\end{itemize} 

In order to understand its definition, let us explain how one can
``see'' characteristics of arbitrary period by using $1$-periodic
orbits of special Hamiltonians. There is a symplectic diffeomorphism
onto a neighbourhood $\mc V$ of the boundary  
$$\Psi : \big( \partial M \times [1-\delta, \, 1], \, d(S\lambda|)
\big)  \longrightarrow ( \mc V, \, \omega), 
\quad \delta > 0 \tx{ small ,}$$
$$\Psi(p, \, S) = \varphi_{\ln (S)} (p) \ ,$$
where $\lambda|$ denotes the restriction of $\lambda$ to $\partial
M$ (we actually have $\Psi ^* \lambda = S  \lambda|$). We define the {\it symplectic
completion} 
$$(\widehat{M}, \, \widehat \omega)  = (M, \, \omega)  \cup_{\,\Psi} \big( \partial M \times
[1, \, \infty[  , \, d(S\lambda|) \big) \ .$$
Consider now Hamiltonians $H$ such that  $H(p, \, S) = h(S)$ for $S \ge
1-\delta$, where $h: [1-\delta, \, 1] \longrightarrow \RR$ is 
smooth. It is straightforward to see that
$$X_H(p, \, S) = -h'(S) X_{\tx{Reeb}}, \quad S \ge 1-\delta \ .$$
The $1$-periodic orbits of $X_H$ that are located on the level $S$
correspond 
to characteristics on $\partial M$ having period $h'(S)$
under the parameterization given by $-X_{\tx{Reeb}}$. The general
principle that can be extracted out of this computation is that {\it one
``sees'' more and more characteristics as the variation of $h$ is
bigger and bigger}.

We define a Hamiltonian to be {\it admissible} if it satisfies $H
\le  0$ on $M$ and it is of the form  $H(p, \, S) =
h(S)$  for $S$ big enough\footnote{Our 
  definition is slightly more general
  than the one in \cite{functors1} in that we prescribe the behaviour
  of admissible Hamiltonians only near infinity and not on the
  whole of $S\ge 1$. It is clear that the a priori $C^0$ bounds 
  require no additional  argument than the one in
  \cite{functors1}. The situation is of course different if one
  enlarges further the admissible class.}, with $h$ convex increasing and such
that there exists $S_0 \ge 1$ with 
$h'$ constant for $S\ge S_0$. We call such a Hamiltonian
  {\it linear at infinity}. Moreover, we assume that the
slope at infinity of $h$ is not the area of a closed characteristic on
$\partial M$ and that all $1$-periodic orbits of $H$ are
nondegenerate. One method to obtain such Hamiltonians is to slightly
perturb functions $h(S)$, where $h:[1-\delta, \,
\infty[ \longrightarrow \RR$ is equal
to zero in a neighbourhood of $1-\delta$, strictly convex on $\{ h > 0
\}$ and linear at infinity with
slope different from the area of any closed characteristic, by a
perturbation localized around the periodic orbits. 
We point out that there are admissible Hamiltonians having arbitrarily large
values of the slope at infinity.

The {\it admissible almost complex structures} are defined to be those
which satisfy the following conditions for {\it large enough values of
  $S$}: 
$$\left\{ \begin{array}{l} 
J_{(p, \, S)}|_\xi = J_0, \\ 
J_{(p, \, S)} \big( \DP{}{S} \big) = \frac 1 {CS} X_{\tx{Reeb}}(p),
\quad C>0, \\ 
J_{(p, \, S)}(X_{\tx{Reeb}}(p)) = -CS \DP{}{S} \ ,
\end{array}\right. 
$$
where $J_0$ is an almost complex structure compatible with the restriction
of $\omega$ to the contact distribution $\xi=\ker \lambda|$ 
on $\partial M$. These
are precisely the almost complex structures which are invariant under
homotheties $(p, \, S) \longmapsto (p, \, aS)$, $a>0$ for large enough
values of $S$. 

The crucial fact is that the function $(p, \, S) \longmapsto S$ is
plurisubharmonic with respect to this class of almost complex
structures. This means that $d(dS \circ J)(v, \, Jv) < 0$ for any
nonzero $v\in T_{(p, \, S)}\big( \partial M \times [1, \, \infty[\big)$, $p\in
\partial M$, $S \ge 1 $ big enough, and indeed we have $d(dS\circ J) =
d(-CS\lambda|) = -C\widehat \omega$. Plurisubharmonicity implies that for any
$J$-holomorphic curve $u:D^2 \longrightarrow \partial M \times [S_0, \,
  \infty [$, $S_0\ge1 $ big enough  one has 
$\Delta(S\circ u) \ge 0$. In particular 
the maximum of $u$ is achieved on the boundary $\partial D^2$ 
(see for example \cite{GT}, Theorem 3.1). A similar argument
applies to solutions of the Floer equation (\ref{Fl eq}), 
as well as to solutions of the parameterized Floer equation (\ref{par
  Fl eq}) for {\it increasing homotopies} satisfying
$\DP{^2h}{s\partial S}\ge 0$~\cite{my survey}. This implies 
that solutions are
contained in an a priori determined compact set. All compactness
results in Floer's theory therefore carry over to this situation  and so does
the construction outlined for closed manifolds. 

%Let us point out that admissible Hamiltonians are time independent by
%definition. The $1$-periodic orbits that correspond to closed
%characteristics can therefore be at most transversally
%nondegenerate. One should perturb the Hamiltonian
%there in order to achieve nondegeneracy. 
%Nevertheless, the pseudoconvexity argument above still applies
%as there are no $1$-periodic orbits at infinity and there is no need
%to perturb $H$ in that region. 

We introduce a partial
order on regular pairs $(H, \, J)$ as 
$$(H, \, J) \prec (K, \, \widetilde{J} ) \qquad \tx{iff} \qquad H\le K
\ . $$
The continuation morphisms (\ref{cont morph}) form a direct system
with respect to this order 
and we define the Floer homology groups as 
$$FH_*(M) = \displaystyle \lim_{\scriptsize 
\begin{array}{c} \rightarrow \\ (H, \,
    J) \end{array} } FH_*(H, \, J) \ .$$
An important refinement of the definition consists in using a
truncation by the values of the action. The latter is decreasing along
Floer trajectories and one builds a $1$-parameter family of
subcomplexes of $FC^*(H, \, J)$, defined as 
$$FC_k^{]-\infty,\, a[}(H, \, J) = \bigoplus_{\scriptsize 
\begin{array}{c} \gamma \in \mc P (H) \\
    i_{CZ}(\gamma) = -k \mod \ 2\nu  \\ A_H(\gamma) < a 
\end{array} } \ZZ \langle \gamma \rangle \ . $$
This allows one to define the corresponding quotient complexes
$$FC_*^{[a, \, b[ } (H, \, J) = FC_*^{]-\infty,\, b[}(H, \, J) /FC_*
^{]-\infty,\, a[} (H, \, J), \quad -\infty \le a < b \le \infty $$
and the same direct limit process goes through. We therefore put
$$FH_*^{[a, \, b[}(M) = \displaystyle \lim_{\scriptsize 
\begin{array}{c} \rightarrow \\ (H, \,
    J) \end{array} } FH_*^{[a, \, b[}(H, \, J) \ .$$

Let us now make a few remarks on the properties of the above invariants. 
\renewcommand{\theenumi}{\alph{enumi}}
\begin{enumerate} 
 \item there is a natural cofinal family of Hamiltonians
 whose values of the action on the $1$-periodic orbits is positive or
arbitrarily close to $0$ (see the construction of $H_1$ in
 Figure \ref{HamPlateau et Rho} (1) described in Section~\ref{la
   preuve de Kunneth}).  
This implies that $FH_*^{[a, \, b[}(M) =0$
 if $b<0$ and $FH_*^{[a, \, b[} $ does not depend on $a$ if the latter
 is strictly negative. In particular we have 
     $$FH_*(M) = FH_*^{[a, \, \infty[}(M), \qquad a<0 \ .$$

\item the infimum $T_0$ of the areas of closed characteristics on the
  boundary is always strictly positive and therefore 
  $$ FH_k^{[a, \, \epsilon[}(M) \simeq \bigoplus _{l \, \equiv \, k
  \mod \ 2\nu }H_{l+m}(M, \, \partial M), $$
  where $a<0\le \epsilon <T_0$, $k \in \ZZ/2\nu \ZZ$, $m = \frac 1 2
  \dim M$. 
  This follows from the fact that, in the limit, the Hamiltonians
  become $C^2$-small on $M \setminus \partial M$ and the Floer complex
  reduces to a Morse complex that computes the relative cohomology, as
  $-\nabla H$ points inward along $\partial M$. 

\item there are obvious truncation morphisms 
  $$FH_*^{[a, \, b[}(H, \, J) \longrightarrow FH_*^{[a', \, b'[}(H, \,
  J), \ a\le a', \ b\le b'$$ 
  which induce morphisms $FH_*^{[a, \, b[}(M) \longrightarrow
  FH_*^{[a', \, b'[}(M)$, $a\le a'$, $b\le b'$.
  If $a=a'<0$, $0 \le b < T_0$ and $b'=\infty$ we obtain a
  natural morphism 
 \begin{equation*} 
 \bigoplus _{l \, \equiv \, k
 \mod \ 2\nu }H_{l+m}(M, \, \partial M) \stackrel{c_*}\longrightarrow FH_k(M) \ ,
 \end{equation*} 
 or, written differently, 
  \begin{equation} \label{c *}
   H_*(M, \, \partial M) \stackrel{c_*}{\longrightarrow} FH_*(M) \ .
  \end{equation} 
  We also note at this point that we have 
  $$FH_*(M) = \displaystyle \lim _{\scriptsize \begin{array}{c}
  \rightarrow \\ b   \end{array} } \lim_{\scriptsize 
\begin{array}{c} \rightarrow \\ (H, \,
    J) \end{array} } FH_*^{[a, \, b[}(H, \, J), \qquad a<0 $$
  and that the two limits above can be interchanged by general
  properties of bi-directed systems. 

\end{enumerate} 

\noindent 
{\bf Fundamental principle (Viterbo~\cite{functors1}).} {\it If the morphism\break 
$H_{*}(M, \, \partial M) \stackrel {c_*}
\longrightarrow FH_*(M) $ is not bijective, then
there is a closed characteristic on $\partial M$. } Indeed, either
there is some extra generator in $FH_*(M)$, or some
Morse homological generator of $H_*(M, \, \partial M)$  is killed in
$FH_*(M)$. 
 The ``undesired guest'' in the first case or the ``killer'' in
the second  case  necessarily corresponds 
to a closed characteristic on $\partial M$. 

%The property of satisfying the SAWC can be relaxed in the following
%way. 

%\medskip 

%\begin{defi}[Viterbo \cite{functors1}] A manifold with contact type
%  boundary satisfies the {\rm Algebraic Weinstein Conjecture (AWC)} if
%  there is a ring of coefficients such that the map $FH^*(M) \stackrel
%  {c^*} \longrightarrow H^{*}(M, \, \partial M)$ is not an
%  isomorphism. 

%  We say that $M$ satisfies case a) of the AWC or case b) of the AWC 
%if $c^*$ is respectively not injective or not surjective. 
%\end{defi} 

%\medskip 

%All three explicit examples given in the introduction are
%interesting in this context. Cotangent bundles often satisfy case a)
%of the AWC but {\it never} satisfy case b), and this last fact 
%gives the main obstruction to exact
%Lagrange embeddings in manifolds satisfying the SAWC. Stein manifolds and negative
%line bundles satisfy the SAWC, but while Stein manifolds are exact, negative
%line bundles are not as the zero section is a closed symplectic
%submanifold. 
    
%\medskip 
%
%The utility of the definition is illustrated by the following result. 
%
%\begin{thm}[Viterbo \cite{functors1}, Thm. 4.1.] \label{appli char}
%  Let $M^{2n}$ be manifold with contact type boundary which satisfies AWCb)
%  in maximal degree i.e. the morphism $FH^n(M) \stackrel {c^*}
%  \longrightarrow H^{2n}(M, \, \partial M)$ is not surjective. Then
%  the same holds for any codimension $0$ submanifold of restricted
%  contact type in $M$. In particular, any closed hypersurface of
%  restricted contact type has a closed characteristic.  
%\end{thm}
%
%\medskip 
%
%
 
\medskip 

The version of Floer homology that we defined above has various
invariance properties \cite{functors1}. The main one that we shall use
is the following. 

\begin{prop} \label{inv Fl h}
 The Floer homology groups $FH_*(M)$ are an invariant of the
 completion $\widehat M$
 in the following sense: for any open set with smooth boundary $U
 \subset M$ such that $\partial U \subset \partial M \times [1, \,
 \infty[$ and the Liouville vector field $S\DP{}{S}$ is transverse and
 outward pointing along  $\partial U$, we have 
 $$FH_*(\widehat M) \simeq FH_*(\widehat U) \ .$$
\end{prop} 

 \demo One can realise a differentiable isotopy between $M$ and $U$
 along the Liouville vector field. This corresponds to an isotopy of
 symplectic forms on $M$ starting from the initial form $\omega =
 \omega_0$ and ending with the one induced from $U$, denoted by
 $\omega_1$. During the isotopy the boundary $\partial M$ remains of
 contact type and the symplectic asphericity condition is preserved. 
An invariance theorem of Viterbo \cite{functors1}
 shows that $FH_*(M, \, \omega) \simeq FH_*(M, \, \omega_1)$. On the other hand  
 $FH_*(M, \, \omega_1) \simeq FH_*(U,
 \, \omega)$ because $(M,\, \omega_1)$ and $(U,\, \omega)$ are symplectomorphic. 
\hfill{$\square$}

\section{Proof of Theorem A}   
\label{la preuve de Kunneth}   
  
Before beginning the proof, let us note that the natural class of
manifolds for which one can {\it define} Floer homology groups for a product
is that of manifolds with restricted contact type boundary. 
The reason is that $\partial \big( M \times N
\big)$  involves the
  full manifolds $M$ and $N$, not only some neighbourhoods of their
  boundaries.  If $X$ and $Y$ are the conformal vector fields on $M$
  and $N$ respectively and $\pi_M: M \times N \longrightarrow M$,
$\pi_N: M \times N \longrightarrow N$  are the canonical projections, 
the natural conformal vector field on 
  $M \times N$ is $Z = \pi ^*_M X + \pi ^*_N
  Y$. 
In order for $Z$ to be defined in a neighbourhood of 
  $\partial ( M \times N)$ 
  it is necessary that $X$ and $Y$ be globally defined.

We only prove a) because 
the proof of b) is entirely dual. One has to reverse arrows and
replace direct limits with inverse limits. The difference in the
statement is due to the fact that the inverse limit functor is 
in general not exact, except when each term of the directed system is a
finite dimensional vector space~\cite{ES}. For the sake of clarity
we shall give the proof under the assumption $\langle c_1(TM), \,
\pi_2(M) \rangle =0$, $\langle c_1(TN), \,
\pi_2(N) \rangle =0$, so that the grading on Floer homology is defined
over $\ZZ$. 

%The crux of the matter is to establish the  
%The existence
%of the morphism from the classical K\"unneth exact sequence
%to~(\ref{suite exacte Kunneth Floer}) will follow 
%easily by truncating the range of the action. 
   
%   Let us bear in mind the Morse homological situation: for a product
%   $M\times N$ of closed manifolds one can consider the Morse function 
%   given by the sum $f(x)+g(y)$ of Morse functions $f:M\longrightarrow
%   \RR$ and
%   $g:N \longrightarrow \RR$. If one
%   chooses on $M \times N$ the product Riemannian metric,  the
%   Morse complex of $f+g$ is the tensor product of the Morse complexes
%   of $f$ and $g$. The algebraic K\"unneth theorem will then give
%   an exact sequence having as middle term the Morse/singular homology
%   of $M\times N$ and as side terms some suitable tensor product of 
%   the homologies of $M$ and $N$. 

\medskip

I. We establish the short exact
sequence~(\ref{suite exacte Kunneth Floer}).
   Here is the sketch of the proof. We consider on $\widehat M \times
   \widehat N$ a Hamiltonian of the form $H(t, \, x) + K(t, \, y)$,
   $t\in \Ss^1$, $x\in \widehat M$, $y\in \widehat N$. For an almost
   complex structure on $\widehat M \times \widehat N$ of the form
   $J^1 \oplus J^2$, with $J^1$, $J^2$ (generic) almost complex
   structures on $\widehat M$ and $\widehat N$ respectively, the Floer
   complex for $H+K$ can be identified, modulo truncation by the
   action issues, with the tensor product of the
   Floer complexes of $(H, \, J^1)$ and $(K, \, J^2)$. Nevertheless,
   the Hamiltonian $H+K$ is not linear at infinity and hence not
   admissible. We refer to~\cite{teza mea} for a discussion of the
   weaker notion of asymptotic linearity and a proof of the fact that $H+K$
   does not even belong to this extended admissible class.   
   The main idea of our proof is to construct an admissible
   pair $(L, \, J)$ whose Floer complex is roughly the same as the one of $(H+K, \,
   J^1 \oplus J^2)$. The Hamiltonian $L$ will have lots of additional
   $1$-periodic orbits compared to $H+K$, 
   but all these will have negative enough action for
   them not to be counted in the relevant truncated Floer
   complexes. 
%Many of these 
%   newly created $1$-periodic orbits will be (very)
%   degenerate but this will have no drawback because
%   they do not enter anyway in the Floer complex.  
 
\medskip 
  
  We define the {\it period spectrum} $\mc S(\Sigma)$ of a contact type hypersurface
  $\Sigma$ in a symplectic manifold as being the set of
  periods of closed characteristics on $\Sigma$, the latter being parameterized by
  the Reeb flow. We assume from now on that the period
  spectra of $\partial M$ and $\partial N$ 
  are {\it discrete and  injective} i.e. the periods of the closed
  characteristics form a strictly increasing sequence, every
  period being associated to a unique characteristic which is
  transversally nondegenerate. This
  property is $C^\infty$-generic among hypersurfaces \cite{T}, while
  Floer homology does not change under a small $C^\infty$-perturbation of
  the boundary (Proposition \ref{inv Fl h}).  
  Assuming a discrete and injective period spectrum amounts therefore to no loss of
  generality.  

  We shall construct cofinal familes of Hamiltonians and almost
  complex structures $(H_\nu, \, J_\nu^1)$, $(K_\nu, \, J_\nu^2)$, $(L_\nu, \,
  J_\nu)$ on $\widehat M$, $\widehat N$ and $\widehat M \times \widehat
  N$ respectively, with the following property. 

\medskip 

{\noindent \bf Main Property.}
{\it Let $\delta >0$ be fixed. For any $b>0$, there is a positive
  integer $\nu(b, \, \delta)$
  such that, for all $\nu\ge \nu(b, \, \delta)$, the following inclusion of
  differential complexes holds:}
\begin{equation} \label{inclusion fondamentale de complexes}   
{\scriptsize
\xymatrix{
\bigoplus_{r+s=k} FC_r^{[-\delta,\, \frac{b}2[}(H_\nu,\, J_\nu^1) \otimes   
      FC_s^{[-\delta,\, \frac{b}2[}(K_\nu,\, J_\nu^2)  \ar[r] & 
      FC_k^{[-\delta,\, b[}(L_\nu,\, J_\nu)    \ar[dl] \\
     \bigoplus_{r+s=k}    
      FC_r^{[-\delta,\, 2b[}(H_\nu,\, J_\nu^1) \otimes   
      FC_s^{[-\delta,\, 2b[}(K_\nu,\, J_\nu^2) \ar[r] &     
      FC_k^{[-\delta,\, 4b[}(L_\nu,\, J_\nu)
}
}
\end{equation}

It is important to note 
that we require the two composed arrows 
{\scriptsize
$$FC_k^{[-\delta,\, b[}(L_\nu,\, J_\nu)
\hookrightarrow FC_k^{[-\delta,\, 4b[}(L_\nu,\, J_\nu),
$$
\begin{equation*} 
\displaystyle{\bigoplus_{r+s=k} FC_r^{[-\delta,\, \frac{b}2[}(H_\nu,\, 
  J_\nu^1) \otimes    
      FC_s^{[-\delta,\, \frac{b}2[}(K_\nu,\, J_\nu^2) }  
\hookrightarrow 
\displaystyle{\bigoplus_{r+s=k}    
      FC_r^{[-\delta,\, 2b[}(H_\nu,\, J_\nu^1) \otimes   
      FC_s^{[-\delta,\, 2b[}(K_\nu,\, J_\nu^2)} 
\end{equation*} 
}

\noindent to be the usual inclusions corresponding to the truncation by the
action. In practice we shall construct {\it autonomous} Hamiltonians
having transversally nondegenerate $1$-periodic orbits, but one
should think in fact of small local perturbations of these ones, 
along the technique of
\cite{CFHW}. The latter consists in perturbing an autonomous
Hamiltonian in the neighbourhood of a transversally
nondegenerate $1$-periodic orbit $\gamma$, replacing $\gamma$ by
precisely  two
nondegenerate $1$-periodic 
orbits corresponding to the two critical points of a Morse function on
the embedded circle given by $\tx{im}(\gamma)$. The Conley-Zehnder
indices of the perturbed orbits differ by one. Moreover, the perturbation
can be chosen arbitrarily small in any $C^k$-norm and the actions of
the perturbed orbits can be brought arbitrarily close to the action of
$\gamma$. 
%The reader may also consult \cite{my survey} for an explicit
%computation of the Conley-Zehnder indices for the orbits
%obtained by perturbing in an analogous manner a transversally
%nondegenerate {\it sphere} of $1$-periodic orbits in $\CC^n$. 

\medskip 

 a. \ Let $S'$, $S''$ be the vertical coordinates on $\widehat M$ and
 $\widehat N$ respectively. Let $(H_\nu)$, $(K_\nu)$ be {\it cofinal} families of
 autonomous Hamiltonians on 
 $\widehat M$ and $\widehat N$, such that $H_\nu(p', \, S')=h_\nu(S')$ for $S'
 \ge 1$, $K_\nu(p'', \, S'') =k_\nu(S'')$ for $S'' \ge 1$, with $h_\nu$, $k_\nu$
 convex and linear of slope $\lambda_\nu$ outside a small neighbourhood
 of $1$. We assumed that the period spectra of $\partial M$ and
 $\partial N$ are discrete and injective and so we can choose  $\lambda_\nu  
   \notin \mc S(\partial M) \ \cup \ \mc S(\partial N)$ 
with $\lambda_\nu \longrightarrow \infty$, $\nu \longrightarrow
 \infty$. We shall drop the subscript $\nu$ in the sequel 
 by referring to $H_\nu$, $K_\nu$ and $\lambda_\nu$ as $H$, $K$ and $\lambda$. 
%   As pointed out before, the sum $H+K$ is not 
%   linear at infinity. We shall deform it in the
%   complement of a suitably chosen compact set, in order to obtain an
%   admissible Hamiltonian $L$ having, roughly speaking, the same
%   truncated Floer
%   complex as $H+K$. The Hamiltonian $L$ will have lots of additional
%   $1$-periodic orbits compared to $H+K$, 
%   but all these will have negative enough action for
%   them not to be counted in the truncated Floer
%   complex $FC_*^{]-\delta, \, b]}$. Many of these 
%   newly created $1$-periodic orbits will be strongly
%   degenerate but this will have no drawback on the proof because
%   they do not enter anyway in the   Floer
%   complex.  
   Let us denote 
   $$\eta_{\lambda} =   
   \tx{dist} \big( \lambda, \, \mc S(\partial M) \, \cup \,   
   \mc S(\partial N) \big) \, > \, 0 \ , $$  
   $$T_0(\partial M)  =  \min \,  
   \mc S(\partial M)  \ , \ \ T_0(\partial N)  =  \min \,  
   \mc S(\partial N) \ ,$$
   $$T_0 = \min \, \big( T_0(\partial M), \, T_0(\partial N) \big) \,
   > \, 0 \ .$$  

\begin{center}  
\begin{figure}[h]   
         \includegraphics{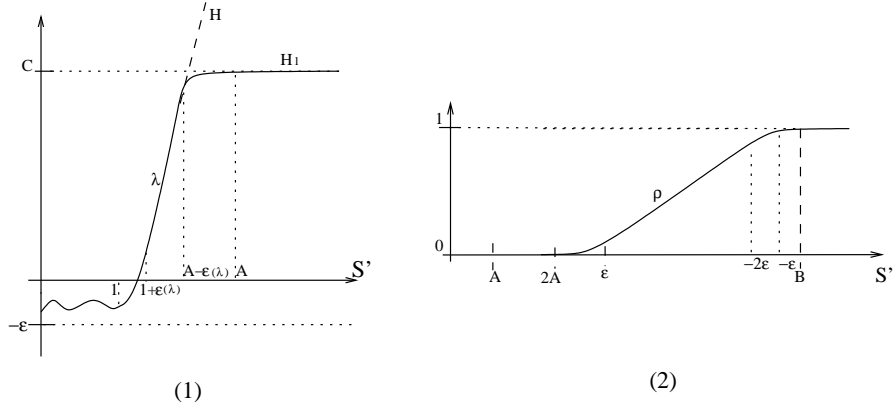}   
  \caption{The Hamiltonian $H_1$ and the truncation function  
   $\rho$ \label{HamPlateau et Rho}}    
\end{figure}    
\end{center}

\medskip 
  
 b. \   Our starting point is the construction by Hermann
   \cite{He} of a cofinal family which allows one to identify, in the
   case of bounded open sets with restricted contact type boundary in
   $\CC^n$, the Floer homologies defined by Viterbo \cite{functors1} and  
   Floer and Hofer \cite{FH}. We fix  
   $$A=A(\lambda)=5\lambda/\eta_\lambda > 1 $$ 
 and consider
   the Hamiltonian $H_1$ equal
   to $H$ for $S' \le  
   A-\epsilon(\lambda)$ and
   constant equal to $C$ for $S' \ge A $, with $C$ arbitrarily close
   to $\lambda(A - 1)$. Here $\epsilon(\lambda)$ is chosen to be small
   enough and positive. We perform the same
   construction in order to get a Hamiltonian $K_1$. 
%The Hamiltonian $H_1$ indeed has additional
%   $1$-periodic orbits compared to  $H$ but, as we shall see below, 
%   these will have very negative action if $A$ is
%   large enough with respect to $\lambda$. 
%
%One suitable choice for $A$ is  
%   \begin{equation}   
%    A = 5\lambda/\eta_\lambda \ .  
%   \end{equation}   
   We suppose (Figure \ref{HamPlateau et Rho} (1))   
   that $H_1$ takes its values in the interval $[-\epsilon,\, 0)$   
   on the interior of $M$ where it is also $C^2$-small 
   and that  $H_1 (\underline{x}, \, S')  
   =h(S')$ on
   $\partial M \times [1,\, \infty[$,  with $h'\equiv\lambda$ on $[1+  
   \epsilon(\lambda),\,  A - \epsilon(\lambda)]$ and $h'\equiv 0$ on  
   $[A,\, \infty[$, where $\epsilon(\lambda)=\epsilon/\lambda$. Thus  
   $H_1$ takes values  
   in $[-\epsilon,\, \epsilon]$ for $S' \in  
   [1,\, 1+\epsilon(\lambda)]$ and in   
   $[\lambda(A-1)-2\epsilon,\, \lambda(A-1)]$ for $S'  
   \in [A-\epsilon(\lambda),\, A]$.   
   
   The Hamiltonian $H_1$ has additional $1$-periodic orbits compared
   to $H$. These are
 either constants on levels $H_1=C$ with action $-C \simeq
 -\lambda(A-1)$, or orbits corresponding to characteristics on
 the boundary, appearing on levels 
  $S=\tx{ct.}$ close to $A$. The action of the latter is arbitrarily
  close to    
    $h'(S)S-h(S) \le (\lambda - \eta_\lambda) \cdot A - \lambda( A -1)  
    +2\epsilon \le -3\lambda \longrightarrow -\infty$, $\lambda  
    \longrightarrow \infty$. The special choice of $A$ is motivated by
    the previous computation. We see in particular that it is crucial to take into
 account the gap $\eta_\lambda$ in order to be able to make the action
 of this kind of orbits tend to $-\infty$.   

c. \ We deform now $H_1 + K_1$  
   to a Hamiltonian that is constant equal to $2C$ 
   outside the compact set $\{ S' \le B, \, S'' \le B \}$, with 
   $$B =   A\sqrt \lambda \ .$$  
   This already holds in $\{ S'  
   \ge A, \, S'' \ge A \}$. We describe the corresponding
   deformation in $\{ S' \le A , \, S'' \ge A \}$ and perform the symmetric
   construction in $\{ S' \ge A, \, S'' \le A \}$. 
   Let us define
   $$H_2 : \widehat{M} \times \partial N \times [A, \, \infty[  
   \longrightarrow \RR \ ,$$  
   $$H_2(\underline{x},\, y,\, S'') = \big( 1-\rho(S'')  \big)  
   H_1(\underline{x}) + \rho (S'') C \ ,$$  
   with $\rho : [A, \ +\infty[ \longrightarrow [0,1]$, $\rho \equiv  
   0 $ on  $[A,\, 2A]$, $\rho \equiv 1$ for $S'' \ge B-\epsilon$,   
   $\rho$ strictly increasing on $[2A,\, B-\epsilon]$,  
   $\rho' \equiv \tx{ct.} \ \in \ [\frac 1 {B-2A-\epsilon}, \frac  
   1 {B - 2A - 3\epsilon} ] $ on $[2A + \epsilon, \, B -  
   2\epsilon]$ (Figure \ref{HamPlateau et Rho} (2)).
   The symplectic form on  $\widehat{M} \times \partial N \times [A,  
   \, \infty[$ is $\omega' \oplus d(S''\lambda'')$, with $\lambda''$  
   the contact form on  $\partial N$. We get   
   $$X_{H_2} (\underline{x},\, y,\, S'') = \big( 1 - \rho(S'') \big)  
   X_{H_1}(\underline{x})  -  
   \big( C - {H_1}(\underline{x}) \big) \rho'(S'') X''_{\tx{Reeb}}(y) \ .$$  
\begin{figure}[h]   
        \begin{center}  
\input{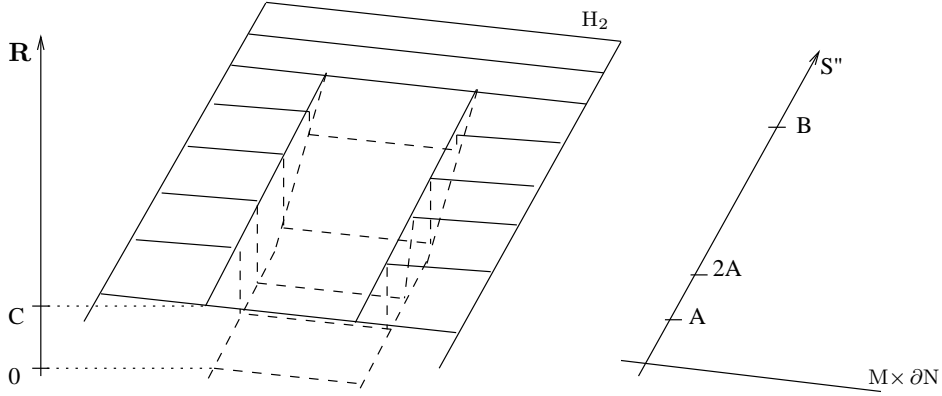} 
\caption{Graph of the deformation $H_2$ \label{deformation}}   
        \end{center}  
\end{figure}    
   The projection of a periodic orbit of  $X_{H_2}$   
   on  $\widehat{M}$  
   is a periodic orbit of $X_{H_1}$. In particular,   
   ${H_1}$ is constant along the projection. Moreover,   the orbits appear on levels  
   $S''=\tx{ct.}$ 
   as there is no component  $\DP{}{S''}$ in  
   $X_{H_2}$. As a consequence, the coefficients in front of
   $X_{H_1}$ and $X''_{\tx{Reeb}}$ are constant along one orbit of  
   $X_{H_2}$.   
   A $1$-periodic orbit $\theta$  of $X_{H_2}$ corresponds therefore  
   to a  couple $(\Gamma, \, \gamma)$ such that   
   \begin{itemize}   
   \item $\Gamma$ is an orbit of  $X_{H_1}$ having period 
   $1-\rho(S'')$;  
   \item $\gamma$ is a closed characteristic on the level 
   $\partial N \times \{ S''\}$,   
   having period  $\big( C - {H_1}(\underline{x})  
   \big) \rho'(S'')$ and having the {\it opposite} orientation than
   the one given by  $X''_{\tx{Reeb}}$. We have used the notation   
   $\underline{x}=\Gamma(0)$.    
   \end{itemize}   
  
   The action of $\theta$ is
   \begin{eqnarray} \label{action totale}   
    A_{H_2+K_1}(\theta) & =  & -A(\Gamma) - A(\gamma) -H_2 - K_1
    \nonumber \\   
    & = & A_{H_1}(\Gamma) -A(\gamma) - \rho (S'') \big( C -  
    {H_1}(\underline{x}) \big) - C   \ .
   \end{eqnarray}   
  
   We have denoted by  $A(\gamma)$, $A(\Gamma)$ the areas of the   
   orbits $\gamma$ and  $\Gamma$ respectively. The Hamiltonian  $K_1$ is  
   constant in the relevant domain and we have directly replaced it by
   its value $C$. It is useful to notice that    
   $A(\gamma)=-S''\rho'(S'')(C-{H_1}(\underline{x}))$. The  
   minus sign comes from the fact that the running orientation on
   $\gamma$ is opposite to the one given by $X''_{\tx{Reeb}}$. 
We have used the symplectic form  $d(S''\lambda'')$ on the
   second factor in~(\ref{action totale}).  

   For $0 \le T = 1-\rho(S'') \le 1$,  
the  $T$-periodic orbits of   
   $X_{H_1}$ belong to one of the following classes.  
   \renewcommand{\theenumi}{\arabic{enumi}} 
   \begin{enumerate}   
    \item \label{11} constants in the interior of  $M$, having zero area;  
    \item \label{22}  
      closed characteristics located around the level $S'=1$, whose
      areas belong to the interval 
      $[T_0(\partial M), \, T\lambda]$;
    \item \label{33}  
      if $T\lambda \in \mc S(\partial M)$, one has a closed
      characteristic of area $S'T\lambda$ for any  $S'  
      \in [1+\epsilon(\lambda),\,  A -\epsilon(\lambda)]$   
      - the interval where  ${H_1}$ is   
      linear of slope  $\lambda$;
    \item \label{44}  
      closed characteristics located around the level $S'=A$, whose
      areas belong to the interval $[T_0(\partial M),\, T(\lambda - \eta_\lambda)A]$;
    \item \label{55} constants on levels $S' \ge A$, having zero area.  
   \end{enumerate}   
 
 %  A suitable choice for the parameter $B$ is 
%   \begin{equation}   
%     B = A \sqrt{\lambda} 
%   \end{equation}   
%   and, 
For $\epsilon >0$ fixed we choose the
   various parameters involved in our constructions such that
   $$\lambda (A-1) \, \ge \, C \, \ge \, \lambda (A-1) - \epsilon, $$
   $$\rho' \, \le  \, 1/ A(\sqrt{\lambda} -1) + \epsilon,$$   
   $$ S'' \tx{ close to } B \ \Longrightarrow \    
   \rho'(S'')\cdot S'' \, \le \, \sqrt{\lambda} / (\sqrt{\lambda}  
   -1) \ .$$   

   We now show that the actions of $1$-periodic orbits of $X_{H_2}$
   appearing in the region   
   $\widehat{M} \times  
   \partial N \times [A,\, \infty[$ tend uniformly to   
   $-\infty$ when $\lambda \rightarrow +\infty$.   
   We estimate the action of the orbits  $\theta$ according to
   the type of their first component  $\Gamma$ and according to 
the level $S''\ge A$ on which lies $\gamma$. 
%   All limits below are taken with respect to
%   $\lambda \longrightarrow +\infty$. 
   \renewcommand{\labelenumii}{\theenumii)}  
   \begin{enumerate}   
   \item  $\Gamma$ of type \ref{11} corresponds to $ {H_1} \in  
     [-\epsilon,0]$. 
    \begin{enumerate}  
     \item \label{premier}   
       $S'' \in [A,2A] \ \bigcup \ [B-\epsilon(\lambda),\,  
       \infty[$. Because $\rho'=0$ there is no component of
       $X_{H_2}$ in the   
       $X''_{\tx{Reeb}}$ direction, orbits  $\Gamma$  
       appear in degenerate families (of dimension $\dim\,N$) 
       and the action of  $\theta$ is   
       $A_{H_2+K}(\theta) \le \epsilon -C $.   
     \item \label{deuxieme}   
       $S'' \in [2A, \, \frac {A + B} 2 ]$. Orbits   
       $\Gamma$ come in pairs with closed characteristics   
       $\gamma $ of period
       $\rho' (C - {H_1}(\underline{x}))$ on  $\partial N \times  
       \{ S'' \}$. We have
       \begin{eqnarray*}   
\lefteqn{A_{H_2+K}(\theta) \ \le \ \epsilon +  
         S''\rho'(S'')(C+\epsilon) -C} \\  
        & \le & \epsilon + \frac{ A+B} 2 \cdot \frac 1 {B - 2A  
         -3\epsilon} \cdot (C+\epsilon) - C \ \le \ 
 -\frac 1 4 C  \   .  
        \end{eqnarray*}  
       The second inequality is valid for sufficiently large
       $\lambda$, in view of  $B =  
       A\sqrt{\lambda}$ which implies  $(B+A)/(B-2A-3\epsilon)  
       \rightarrow 1$.      
      \item \label{troisieme}   
        $S'' \in [\frac {A + B} 2, B -\epsilon(\lambda)]$ (hence
        $\rho \in [\frac 12, 1]$). For $\lambda $ big enough we have    
        \begin{eqnarray*}   
         \lefteqn{A_{H_2+K}(\theta) \ \le \ \epsilon +  
         S''\rho'(S'')(C+\epsilon) - \rho(S'')C - C} \\  
         & \le & \epsilon + \frac{B - \epsilon}{B - 2A -  
         3\epsilon} \cdot (C+\epsilon) - \frac 12 \cdot C -C \ \le \ 
         -\frac 14 C \ .   
        \end{eqnarray*}  
     \end{enumerate}   
   \item $\Gamma$ of type \ref{22} corresponds to ${H_1} \in [-\epsilon,  
     \epsilon]$ and $S' \in [1, \ 1+\epsilon(\lambda)]$.   
     The area of $\Gamma$ belongs to  the interval 
     $[T_0(\partial M), \, (1-\rho(S'')) \lambda]$. 
      \begin{enumerate}   
        \item $S'' \in [A,\, 2A] \ \bigcup \ [B-\epsilon(\lambda),\,  
       \infty[$. As in \ref{premier}) we have
%As in \ref{premier}) orbits $\theta$ appear in
%       degenerate families and their total action is    
       $$A_{H_2+K}(\theta) \ \le \ \epsilon + (1-\rho(S''))\lambda - C  
       \ \le \ \epsilon + \lambda - C \ .$$  
        \item $S'' \in [2A, \frac {A + B} 2  
          ]$. Like in \ref{deuxieme}) the total action of  $\theta$  
          is
          \begin{equation*}  
           A_{H_2+K}(\theta)  \le   \epsilon +  
           (1-\rho(S''))\lambda + S''\rho'(S'')(C+\epsilon) -C  \le
             - \frac 12 C   \, .  
          \end{equation*}    
        \item $S'' \in [\frac {A + B} 2 , \, B
 -\epsilon(\lambda)]$. Following \ref{troisieme}) one has   
          \begin{equation*}  
           A_{H_2+K}(\theta) \le  \epsilon +(1-\rho(S''))\lambda  
           + S''\rho'(S'')(C+\epsilon) -\rho(S'')C - C \le -\frac 14 C  .   
          \end{equation*}    
      \end{enumerate}   
    \item $\Gamma$ of  type \ref{33} has an  action   
          $A_{H_1}(\Gamma) \le S'T\lambda - \lambda(S'-1-\epsilon')  
          \le (1+\epsilon') \lambda$, where  
          $\epsilon'=\epsilon(\lambda)$.    
          \begin{enumerate}   
            \item $S'' \in [A,\, 2A] \, \bigcup \,  
       [B-\epsilon(\lambda),\, \infty[$~:    
            $A_{H_2+K}(\theta) \le 2\lambda - C$.  
            \item $S'' \in [2A, \, \frac {A + B} 2 ]$~:    
             $A_{H_2+K}(\theta) \le 2\lambda -\frac 14  C$.   
            \item \label{quatrieme} $S'' \in [\frac {A + B} 2 , \, B  
              -\epsilon(\lambda)]$. The technique used in   
              \ref{troisieme}) in order to get the upper bound no longer applies, as   
              $C-{H_1}(\underline{x})$ can be arbitrarily close to
              $0$. Nevertheless $\rho$ satisfies by 
            definition the inequality $(S''-2A) \rho'(S'') \le \rho(S'')+\epsilon$.
              We thus get 
              \begin{eqnarray*}   
              \lefteqn{A_{H_2+K}(\theta)} \\
              \ & \le &  (1+\epsilon')\lambda +  
              S''\rho'(S'')(C-{H_1}(\underline{x})) -  
              \rho(S'')(C-{H_1}(\underline{x})) -C \\   
              \  & \le &  (1+\epsilon')\lambda +  
              \frac{2A}{B-2A-3\epsilon} (C+\epsilon)  
              +\epsilon(C+\epsilon) -C \le 2\lambda - \frac 12 C \, .  
              \end{eqnarray*}    
             \end{enumerate}    
    \item $\Gamma$ of type \ref{44} corresponds to $ {H_1} \in  
      [C-\epsilon,\, C]$ and $A_{H_1}(\Gamma)\le (1-\rho(S''))  
      \lambda A-\lambda(A-1) \le \lambda$. In all three
      cases a)-c) we get $A_{H_2+K}(\theta) \le 2\lambda - \frac 1 4
      C$. 
% in fact, $A_{H_1}(\Gamma)$  
%      goes to $-\infty$, but we don't need such a strong piece of 
%      information).   
%        The total action of  $\theta$ is respectively   
%        
%          \item $A_{H_2+K}(\theta) \le \lambda - \frac 14 C \  
%            \longrightarrow \ -\infty$~;   
%          \item $A_{H_2+K}(\theta) \le 2\lambda - \frac 12 C  
%            \ \longrightarrow \ -\infty$.    
%        \end{enumerate}   
    \item $\Gamma$ of type  \ref{55} corresponds to ${H_1} \equiv C$. Like
      in   
      \ref{premier}) there is no component in the  $X''_{\tx{Reeb}}$
      direction for $X_{H_2}$  
      and orbits  $\Gamma$ appear in (highly) degenerated families. The
      total action in all three cases a) - c) is   
       $A_{H_2+K}(\theta) = -C - C = -2C$.  
    \end{enumerate}

  This finishes the proof of the fact that the action of the new orbits of $H_2$
  tends uniformly to $-\infty$. 
%The direct
%  consequence is that they are not going to be involved in the relevant Floer
%  complex that computes the (co)homology of  
%  $\widehat M \times \widehat N$. 

\begin{center}  
\begin{figure}[h]   
         \includegraphics{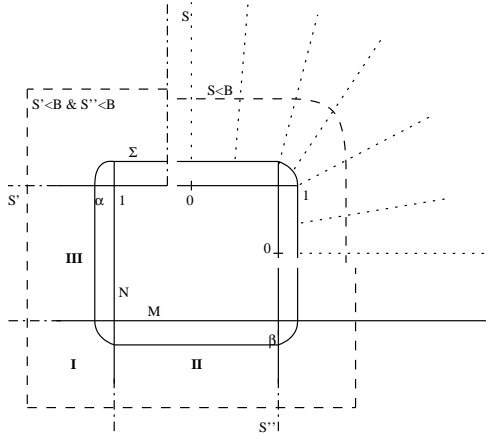}   
  \caption{Parameterization of the product $\widehat M \times \widehat
    N$ \label{param product}}    
\end{figure}    
\end{center}

\medskip 

  d. \ The symmetric construction can be carried out for  $K$ in the
    region   
    $\partial M \times [A, \, \infty[ \times \widehat{N}$. One gets
    in the end  a Hamiltonian   
    $H_2+K_2$ which is constant equal to  $2C$ on $\{  
    S' \ge B \} \bigcup \{ S'' \ge B \}$. We  modify  
    now $H_2+K_2$ outside the compact set  $\{S'\le B \} \bigcap \{  
    S'' \le  B \}$ in order to make it linear with respect to the
    Liouville vector field  $Z =
    X \oplus   
    Y$ on $\widehat{M} \times \widehat{N}$.   
    Let us define the following domains in $\widehat M \times \widehat
    N$ (see Figure \ref{param product}): 
$$\tx{\bf I} = \partial  M \times   
 [1,\ +\infty[ \  \times \ \partial N \times [1, \ +\infty[ \ ,$$   
 $$\tx{\bf II} = M \times \partial N \times [1, \ +\infty[ \ , \ \ 
\tx{\bf III} = \partial M \times [1, \ +\infty[ \ \times \  N \ .$$   
Let $\Sigma \subset \widehat{M} \times \widehat{N}$ be a hypersurface 
which is transversal to $Z$ such that 
    $$S'\ _{\vert_{\Sigma \ \cap \ \text{\bf III}}} \equiv \alpha > 1,  
    \qquad \qquad   
    S'\ _{\vert_{\Sigma \ \cap \ \text{\bf I}}} \in [1,\ \alpha] \ ,$$  
    $$S''\ _{\vert_{\Sigma \ \cap \ \text{\bf II}}} \equiv \beta > 1,  
    \qquad \qquad   
    S''\ _{\vert_{\Sigma \ \cap \ \text{\bf I}}} \in [1 , \ \beta] \ .  
    $$  
    We parameterize $\widehat{M} \times \widehat{N} \ \setminus \  
    \tx{int}(\Sigma)$ by   
    $$\Psi~: \Sigma \times [1, \ +\infty[ \longrightarrow \widehat{M}  
    \times \widehat{N} \ \setminus \ \tx{int}(\Sigma) \ , $$  
    $$(z,S) \longmapsto \big( \varphi'_{_{\ln S}}(z), \ \varphi''  
    _{_{\ln S}} (z) \big) \ ,$$  
    which is a  symplectomorphism if one endows $\Sigma \times [1,\  
    +\infty[$ with the symplectic form $d(S \lambda|)$,
    where $\lambda| = \iota(X\oplus Y)\big( \omega \oplus
    \sigma \big)\vert_\Sigma$. As 
    an example, for $z \in \Sigma \ \cap \ \text{\bf III}$ we have   
    $\varphi' _{_{\ln S}}(z) = \big( x(z), \ S\alpha \big)$. It is
    easy to see that    
    \begin{equation} \label{estimation}   
    \Psi ^{-1} \Big( \{  
    S' \ge B \} \ \cup \ \{ S'' \ge B \} \Big) \supseteq \{ S \ge  
    B \} \ .  
    \end{equation}   
    As a consequence $H_2+K_2$ is constant equal to  $2C$ on 
    $\{ S \ge B   
    \}$. We replace it by 
    $L =l(S)$ on $\{ S \ge B\}$, with $l$ convex and   
    $l'(S) = \mu \notin \mc S(\Sigma)$   
    for $S \ge B+\epsilon$. The additional $1$-periodic orbits that are
    created in this way have action $A_L \le \mu(B+\epsilon)  
    -2C = \mu(\sqrt{\lambda}A+\epsilon) - 2\lambda (A-1)$.   
    By choosing $\mu = \sqrt{\lambda}$ one ensures  
    $A_L \longrightarrow \ -\infty$, as well as the cofinality of
    the family of Hamiltonians $L$ as $\lambda\rightarrow \infty$. Indeed,
    the Hamiltonian  $L$ is bigger than  $(\sqrt{\lambda}  
    -\epsilon)(S-1)$ on $\Sigma \times [1, \, \infty[$.   
    Note that the choice of  $\mu$ equal to    
    $\sqrt{\lambda}$ and not belonging to the spectrum of  $\Sigma$  
    is indeed possible if we choose $\Sigma$ to have a discrete and
    injective spectrum.    
%\vspace{-.9cm}
\quad \ e. \ We have constructed a cofinal family of Hamiltonians
  $(L_\nu)_{\nu\ge 1}$ which are linear at infinity and that are
  associated to the initial Hamiltonians $(H_\nu)$ and
  $(K_\nu)$.   
  I claim that the {\it Main Property} holds for $(L_\nu)_{\nu \ge
    1}$. We assume of course that 
 $J_\nu^1$, $J_\nu^2$ and $J_\nu$ 
are regular almost complex structures  for 
$H_\nu$, $K_\nu$, $L_\nu$ respectively, which are standard  
for $ S' \ge 1+\epsilon$, $S'' \ge 1+\epsilon$ and $S\ge  
B+\epsilon$. Moreover, the almost complex structure $J_\nu$ is of the form 
 $J_\nu = J_\nu^1 \oplus J_\nu^2$ for $S \le  
B$. The preceding estimates on the action of the $1$-periodic orbits
  show that  the sequence (\ref{inclusion fondamentale
     de complexes}) of inclusions is certainly valid at the level of
  modules. This says  in particular that the $1$-periodic orbits involved
  in the free modules appearing in (\ref{inclusion fondamentale
     de complexes}) are located in a neighbourhood of 
 $M$, $N$ and  $M\times N$ respectively, and we can assume that the
  latter is contained in  $\{ S \le 1  
\}$ for a suitable choice of  $\Sigma$. 
Classical transversality arguments  (\cite{FH}, Prop. 17~;  
\cite{FHS}, 5.1 and 5.4) ensure that we can choose regular
  almost complex structures of the form described above. 
The point is now to prove that the inclusions~(\ref{inclusion fondamentale  
  de complexes}) are also valid at the level of {\it differential}
  complexes. 
%To this purpose we shall use an energy estimate for Floer
%  trajectories (Lemma~\ref{petit lemme energetique}). 
%
%  \medskip 

It is enough to show that any  trajectory 
 $u = (v, \, w):  
\RR \times \Ss^1 \longrightarrow \widehat{M} \times \widehat{N}$
 which satisfies
 $u(s,\cdot) \longrightarrow (x^{\pm}, \, y^{\pm})$,  
$s\longrightarrow \pm \infty$ stays in the domain  $\{ S\le 1  
\}$, where the Floer equation is split. That will  imply that $v$
 and $w$ are Floer trajectories in 
$\widehat{M}$ and $\widehat{N}$ respectively and will  prove   
(\ref{inclusion fondamentale de complexes}) at the level of differential complexes.  
It is of course enough to prove this statement under the assumption
 that  
$J_\nu=J_\nu^1\oplus J_\nu^2$ on the whole of $\widehat{M} \times
 \widehat{N}$: 
the  Floer trajectories would then be a-posteriori contained in $\{
 S\le 1 \}$. Moreover, the proof of this fact only makes use of the split
 structure on the set $\{ S' \le 2A, \, S'' \le 2A \}$ and this 
 allows one to modify  $J_\nu$ outside $\{ S \le
 B \}$ in order to formally work with an almost complex structure that
 is homothety-invariant at infinity. 
  
We  therefore suppose  in the sequel that  $J_\nu = J_\nu ^1 \oplus  
J_\nu^2$. Arguing by contradiction, assume that the Floer trajectory
$u$ is not contained in 
$\{ S  
\le 1 \}$. As  $u$ is anyway contained in a compact set, we infer that
the function  $S \circ u$ has a  local maximum in  $\{ S  
> 1 \}$, which means that one of the functions  $S'\circ v$ or  $S''
\circ w$ has a local maximum in $\{ S' > 1+\epsilon \}$ or $\{ S''  
> 1+\epsilon \}$ respectively. The two cases are symmetric and we can
assume without loss of generality that 
$S'' \circ w$  
has a local maximum in  $\{ S'' > 1+\epsilon \}$.  
In view of the fact that  $w$ satisfies the Floer equation associated
to  $K_\nu$ for $S'' \le 2A$ 
the maximum principle ensures that the value of $S'' \circ
w$ at the local maximum is in the interval $]2A, \,  
\infty[$.  But  $w(s,\cdot) \rightarrow y^\pm$, $s \rightarrow \pm \infty$  
with $y^\pm \in \{ S'' \le 1+\epsilon  \}$ and this implies  that $w$  
crosses the  hypersurfaces $\{S'' =A \}$ and $\{ S'' =2A\}$. Moreover,   
the piece of  $w$ contained in  $\{ A \le S'' \le 2A \}$ is    
$J_\nu^2$-holomorphic because  $K_\nu$ is constant on that strip.   
We therefore obtain
\begin{eqnarray} \label{inegalite estimation energie}   
\lefteqn{A_{L_\nu}(x^-, \, y^-) -A_{L_\nu}(x^+, \, y^+) \ = \ \int_{\RR \times  
  \Ss^1} |(v_s,\, w_s)|^2_{J_\nu^1 \oplus J_\nu^2}}  \\
& = & \int_{\RR \times\Ss^1} |v_s |^2_{J_\nu^1}+|w_s|^2_{J_\nu^2} \
 \ge \ \int_{\RR \times \Ss^1} |w_s|^2_{J_\nu^2} \nonumber \\
& \ge &   
\int_{\big[(s,\, t) \,~: \, w(s,\, t) \, \in \, \{ A \le S'' \le 2A \}  
  \big]} |w_s|^2_{J_\nu^2} \ = \ \tx{Area}(w \cap \{ A \le S'' \le 2A
\}) \ . \nonumber 
\end{eqnarray}

The last equality holds because $w$ is 
$J_\nu^2$-holomorphic in the relevant region. 
The lemma below will allow us to conclude. 
It is inspired from Hermann's work \cite{He},  
where one can find it stated  for $\widehat{N} = \CC^n$.  

\medskip 
  
\begin{lem} \label{petit lemme energetique}   
Let $(N, \, \omega)$ be a manifold with contact type boundary and 
$\widehat{N}$ its symplectic completion. 
  Let $J$ be an almost complex structure which is homothety invariant
  on  
  $\{ S \ge 1 \}$. 
  There is a constant  $C(J) >0$ such that,  
  for any  $A \ge 1$ and any  $J$-holomorphic curve $u$ having
  its boundary components on both   
  $\partial N  \times \{ A \}$ and $\partial N  \times \{2A \}$, one has    
  \begin{equation} \label{petite ineg aires}  
    \tx{Area}(u) \ge C(J) A \ .  
  \end{equation}   
\end{lem}   
  
\demo 
%The key point of the proof is the fact that a standard almost
%complex structure is invariant under dilatations in the vertical
%coordinate. More precisely, 
Consider, for  $A \ge 1$, the map    
$$\begin{array}{rcl}
\partial N  \times [1, \, \infty[ & \stackrel {h_A}   
\longrightarrow & \partial N \times [1, \,  
\infty[ \ , \\  
& & \\
(p, \, S) & \longmapsto & (p, \, AS) \ .
\end{array}
$$  
By definition the map $h_A$ is $J$-holomorphic. On the other
hand, by using the explicit form of $J$ given in
\S\ref{constructions}, one sees that $h_A$ 
expands  the area
element by a factor  $A$.  As a consequence, up to rescaling by $h_A$
it is enough to prove~(\ref{petite  
  ineg aires})  for  $A=1$. We apply Gromov's
Monotonicity Lemma \cite{G} 1.5.B, 
%\cite{Muller} 4.2.1, 
\cite{Sik} 4.3.1  
which ensures the existence of $\epsilon_0>0$ and of
$c(\epsilon_0, \, J)>0$ 
such that,  
for any $0 < \epsilon \le \epsilon_0$, any $x\in \partial N  \times  
[1, \, 2]$ with  $B(x, \, \epsilon) \subset \partial N  \times [1, \,  
2]$ and any connected 
$J$-holomorphic curve $S$ such that $x\in S$ and $\partial
S\subset \partial B(x,\, \epsilon)$ one has    
$$\tx{Area}(S \cap B(x, \, \epsilon)) \ge c(\epsilon_0, \, J) 
\epsilon^2 \ .$$  
Let us now fix $\epsilon$ small enough so that  $B(x, \,  
\epsilon) \subset \partial N  \times [1, \, 2]$ for all  $x \in \partial N   
\times \{ \frac 3 2 \}$. As the boundary of $u$ rests on both   $\partial N   
\times \{ 1 \}$ and $\partial N  \times \{ 2 \}$ one can find such a point   
 $x$ on the image of  $u$. We infer  $\tx{Area}(u) \ge  
\tx{Area}(u \cap B(x, \, \epsilon)) \ge c(\epsilon_0, \, J) 
\epsilon^2$. Then $C(J)=c(\epsilon_0, \, 
J)\epsilon^2$ is the desired constant.
      \hfill{$\square$}  
  
\medskip 

Applying  Lemma  \ref{petit lemme energetique} to our situation
we get a constant   
$c >0$ that does not depend on $u$ and such that $\tx{Area}(w \cap  
\{ A \le S'' \le 2A \} ) \ge cA \ge C\lambda$. The difference of the
actions  in (\ref{inegalite estimation energie}) is at the same time
bounded by   
 $4b$ and we get a contradiction for  $\lambda$  
large enough.  
 
For a fixed $b>0$ the Floer trajectories corresponding to  $L_\nu$ are therefore  
contained in  $\{ S \le 1 \}$ for $\nu$ large enough. This
proves that the sequence of inclusions  (\ref{inclusion fondamentale
  de complexes}) is valid at the level of differential complexes.
A few more commutative diagrams will now finish the proof. 
First, as a direct consequence of (\ref{inclusion
  fondamentale de complexes}), one has the commutative diagram   
\begin{equation} \label{premier diagramme intermediaire}    
{\scriptsize
\xymatrix{   
 H_k \big( FC_*^{[-\delta,\, \frac{b}2[}(H_\nu,\, J_\nu ^1) \otimes   
      FC_*^{[-\delta,\, \frac{b}2[}(K_\nu,\, J_\nu ^2)   
\big) \ar@{-->}[d] \ar[r]   & FH_k^{[-\delta,\, b[}(L_\nu,\, J_\nu)  
\ar[dl] \ar@{-->}[d] \\
H_k \big( FC_*^{[-\delta,\, 2b[}(H_\nu,\, J_\nu ^1) \otimes   
      FC_*^{[-\delta,\, 2b[}(K_\nu,\, J_\nu ^2)   
\big) \ar[r] &   FH_k^{[-\delta,\, 4b[}(L_\nu,\, J_\nu)     
}   
}
\end{equation}     
 By taking the direct limit for $\nu \rightarrow
\infty$ and $b \rightarrow \infty$ this induces the diagram

\begin{equation} \label{deuxieme diagramme intermediaire}   
{\scriptsize  
\xymatrix{ \displaystyle   
 \lim_{\stackrel{\longrightarrow} {b \rightarrow +\infty}} 
\lim_{\stackrel{\longrightarrow} {\nu {\longrightarrow} +\infty}}    
H_k \big( FC_*^{[-\delta, \frac{b}2[}(H_\nu, J_\nu ^1) \otimes   
      FC_*^{[-\delta, \frac{b}2[}(K_\nu, J _\nu ^2)   
\big) \ar[r] \ar@{-->}[d]  &   FH_k(M\times N) 
\ar[dl]_\sim \ar@{-->}[d] \\ 
 \displaystyle   
 \lim_{\stackrel{\longrightarrow} {b \rightarrow +\infty}}
 \lim_{\stackrel{\longrightarrow} { \nu {\longrightarrow} +\infty}}   
H_k \big( FC_*^{[-\delta, 2b[}(H_\nu, J_\nu ^1) \otimes   
      FC_*^{[-\delta, 2b[}(K_\nu, J_\nu ^2)   
\big) \ar[r]  &  FH_k(M \times N)     
}  
}
\end{equation}   
 
We easily see that {\it vertical arrows are isomorphisms} as the 
direct limits following $\nu$ and $b$  commute with each
other (this is a general property of bidirected systems). 
This implies that {\it the diagonal arrow is  an isomorphism} as
well. At the same time the algebraic K\"unneth
theorem~\cite[VI.9.13]{D} ensures the existence of a split short exact
sequence 
\begin{equation} \label{troisieme diagramme intermediaire}   
{\scriptsize  
\xymatrix{     
  \bigoplus_{r+s=k} FH_r(M) \otimes    
  FH_s(N) \ \ \ar@{>->}[r] &    
  \displaystyle    
   \lim_{\begin{array}{c} \scriptsize \longrightarrow \\ b,\, \nu
       \rightarrow \infty \end{array}
}
H_k \big( FC_*^{[-\delta,\, 2b[}(H_\nu,\, J _\nu ^1) \otimes   
      FC_*^{[-\delta,\, 2b[}(K_\nu,\, J _\nu ^2)   
\big)  \ar@{->>}[d]  \\    
  & \bigoplus_{r+s=k-1} \tx{Tor}_1^A \big( FH_r(M), \    
  FH_s(N) \big) 
}     
}
\end{equation}     

We infer the validity of the short exact sequence~(\ref{suite exacte
  Kunneth Floer}). 
%{\scriptsize
%\begin{equation}
%\xymatrix{     
%  \bigoplus_{r+s=k} FH_r(M) \otimes    
%  FH_s(N) \ \ \ar@{>->}[r] &    
%  \displaystyle    
%   FH_k(M\times N)   
%  \ar@{->>}[d] \\ 
%  & \bigoplus_{r+s=k-1} \tx{Tor}_1 \big( FH_k(M), \    
%  FH_s(N) \big) 
%}  
%\end{equation}    
%}
We note that we use in a crucial way the
exactness of the direct limit functor in order to 
obtain~(\ref{troisieme diagramme intermediaire}) 
from the K\"unneth exact sequence in truncated homology.

\medskip   
  
II. We prove now the existence of the morphism from the classical
K\"unneth exact sequence to~(\ref{suite exacte Kunneth Floer}). 
We restrict the domain of the action to  $[-\delta, \,
    \delta[$ with $\delta >0$ small enough.
    The Floer trajectories of a 
    $C^2$-small autonomous Hamiltonian
    which is a Morse function on $M$ coincide in the symplectically
    aspherical case with the gradient trajectories in the
    Thom-Smale-Witten complex. We denote by $C_*^{\tx{Morse}}$ the
    Morse complexes on the relevant manifolds. By~(\ref{inclusion
      fondamentale de complexes}) there is a commutative diagram
\begin{equation} \label{premier diag pour Kunneth geant}  
{\scriptsize  
\xymatrix{  
\displaystyle{\bigoplus_{r+s=k} FC_r^{[-\delta, 2b[}(H_\nu,  J_\nu^1) \otimes   
      FC_s^{[-\delta, 2b[}(K_\nu,  J_\nu^2) } \ \ \ar@{^(->}[r]
      & FC_k^{[-\delta, 4b[}(L_\nu,  J_\nu)   \\  
\displaystyle{\bigoplus_{r+s=k} C^{\tx{Morse}}_{m+r}(H_\nu,  
J_\nu^1)   
\otimes C^{\tx{Morse}}_{n+s}(K_\nu, J_\nu^2) }
 \ar[u] 
\ar@{=}[r] & C^{\tx{Morse}}_{m+n+k}(L_\nu, J_\nu^1\oplus J_\nu ^2) \ .
\ar[u]  
}  
}
\end{equation}   
The relevant Morse complexes compute homology relative to the
boundary. With an obvious notation the above diagram induces in homology
\begin{equation} \label{deuxieme diag pour Kunneth geant}  
{\scriptsize  
 \xymatrix{  
\displaystyle \lim_{\stackrel{\longrightarrow } b } \ \lim_{\stackrel 
\longrightarrow n} \  
H_k\big(FC_*^{[-\delta,\, 2b[}(H_\nu, \, J_\nu^1) \otimes   
      FC_*^{[-\delta,\, 2b[}(K_\nu, \, J_\nu^2)\big)    
\ar[r]^{\qquad \qquad \qquad \qquad \qquad \sim}  
      & FH_k(M\times N)   \\  
H_{m+n+k} \big( C_{m+*}(M, \, \partial M) \otimes   
C_{n+*}(N, \partial N) \big) \ar[u]^{\phi}  
\ar[r]^{\qquad \qquad \sim}  
 & H_{m+n+k}(M\times N, \, \partial(M\times N))  \ar[u]^{c_*} 
}  
}
\end{equation}   
By naturality of the algebraic K\"unneth exact sequence, the map $\phi$ 
fits into the diagram below, where $\, * \, $ stands for
its domain and target.
\begin{equation} \label{troisieme diag pour Kunneth geant}  
{\scriptsize  
\xymatrix{   
{\displaystyle{\bigoplus_{r+s = k}}} FH_r(M) \otimes    
  FH_s(N) \ \ \ar@{>->}[r]   &   
    {*} \ar@{->>}[r]   
&  {\displaystyle{\bigoplus_{r+s=k-1}}} \tx{Tor}_1^A \big( FH_r(M),    
  FH_s(N) \big)  \\  
\hspace{-.1cm}{\displaystyle{\hspace{.03cm}\bigoplus_{r+s = k}}}
\hspace{-.1cm}
H_{m+r}(M, \partial M) \otimes    
  H_{n+s}(N, \partial N) \ \ \ar@{>->}[r] \ar[u]_{c_*\otimes c_*} &  
    {*} \ar@{->>}[r]  \ar[u]_{\phi} &   
\hspace{-.32cm}{\displaystyle{\bigoplus _{r+s=k-1}}} \hspace{-.25cm} 
\tx{Tor}_1^A \big(H_{m+r}(M,\partial M),
H_{n+s}(N,\partial N) \big)
\ar[u]_{\tx{Tor}_1(c_*)}  
}  
}
\end{equation}   
Diagrams
(\ref{deuxieme diag pour Kunneth geant}--\ref{troisieme diag pour
  Kunneth geant}) establish the 
desired morphism of exact sequences. 
\hfill{$\square$}

\section{Applications} \label{appli} 
 
%This section contains the proofs of Theorems B and C, as well as some
%other applications of the K\"unneth exact
%sequence related to the existence of closed
%characteristics on hypersurfaces with contact type boundary, explicit
%computations of Floer homology groups and computations 
%of symplectic capacities.

\subsection{Computation of Floer homology groups} 

\begin{prop}  \label{prop:produc with C}
  Let  $N$ be a compact symplectic manifold with restricted contact
  type boundary and let $\widehat N$ be its symplectic completion.  
  Let $\widehat N \times \CC^\ell$, $\ell\ge 1$ be endowed with the
  product symplectic form. Then 
  $$FH_*(\widehat N \times \CC^\ell) =0 \ .$$  
\end{prop}   
  
\demo This follows directly from the K\"unneth exact sequence 
and from Floer, Hofer and Wysocki's computation $FH_*(\CC^{\ell})=0$~\cite{FHW}.
\hfill{$\square$}  

\medskip 

We now fix terminology for 
the proof of Theorem C following~\cite{Eli-psh}. 
A Stein manifold $V$ is a triple $(V,\, J_V,\, \phi_V)$,
where $J_V$ is a complex structure and $\phi_V$ 
is an exhausting plurisubharmonic function. 
We say that $V$ is of {\it finite
type} if we can choose $\phi_V$ with all critical points lying in 
a compact set $K$. The Stein domains $V_c=\{\phi_V\le c\}$ such that
$V_c\supset K$  are called {\it big Stein domains} of $\phi_V$~; they are
all isotopic.
%We say that $V$ is {\it complete} if the gradient of
%$\phi_V$ with respect to the metric $-d(d\phi_V\circ J_V)(\cdot,\,
%J_V\, \cdot)$ is complete in positive time. 
We call $( V,
\, J_V,\, \phi_V)$ {\it subcritical} if $\phi_V$ is Morse and all its
critical points have indices 
strictly smaller than $\frac 1 2 \dim _\RR V$ (they are anyway at most
equal to $\frac 1 2 \dim _\RR V$).   

Let $(V,\, J_V,\, \phi_V)$ be a Stein manifold of finite
type with $\phi_V$ Morse. Following~\cite{SS} we define a 
{\it finite type Stein deformation} of $( V,\, J_V,\, \phi_V)$ 
as a smooth family of complex structures $J_t$,
$t\in [0,1]$ together with exhausting plurisubharmonic functions
$\phi_t$ such that: i) $J_0=J_V$, $\phi_0=\phi_V$; ii) the $\phi_t$ have only Morse or birth-death
type critical points; iii) there exists $c_0$ such that all $c\ge c_0$
are regular values for $\phi_t$, $t\in [0,1]$. 
Condition ii) is not actually imposed in~\cite{SS}, 
but we need it for the following theorem, 
which says that the existence of finite type  
Stein deformations on subcritical
manifolds is a topological problem. 

\medskip

{\bf Theorem (compare~\cite[3.4]{Eli-psh}).} {\it Let $(J_0,\, \phi_0)$ and
$(J_1,\, \phi_1)$ be
finite type Stein structures on $V$.
Assume $J_0$, $J_1$ are homotopic as almost complex
structures and $\phi_0$, $\phi_1$ can be connected 
by a family $\phi_t$, $t\in [0,1]$ of exhausting smooth functions
which satisfy
conditions ii)-iii) above and whose nondegenerate critical points
have subcritical index for all $t\in[0,1]$.
%such that the
%sublevel sets $\{f_0\le c\}$ and   
%$\{f_1\le c\}$ are diffeomorphic for $c$ large enough.
%there exists a family of exhausting smooth
%functions $f_t:V\longrightarrow \RR$, $t\in[0,1]$ such that the above conditions
%ii)-iii) hold and $f_t$ has no critical points of index $n$ for all $t\in[0,1]$.
Then $(J_0,\,
\phi_0)$, $(J_1,\, \phi_1)$ are homotopic by a finite type Stein deformation.}

\medskip

This is the finite type version of Theorem 3.4 in~\cite{Eli-psh}. It
holds because the
latter is proved by \emph{h}-cobordism methods within the
plurisubharmonic category~\cite[Lemma
3.6]{Eli-psh}, and these preserve the
finite type condition. 

\medskip 

Given a Stein manifold of finite type $(V,\, J_V,\, \phi_V)$, let us
fix $c_0$ such that all $c\ge c_0$ are regular values of $\phi_V$. Let
$\omega_V=-d(d\phi_V\circ J_V)$. The
big Stein domains $V_c=\{\phi_V\le c\}$, $c\ge c_0$ are diffeomorphic and
endowed with exact symplectic forms $\omega_c = \omega_V|_{V_c}$
for which $\partial V_c$ is of restricted contact type.
% (the
%Liouville vector field is the dual of $-d\phi_V\circ J_V$ with respect
%to $\omega_c$) 
The Floer homology groups $FH_*(V_c)$ are well defined
and, by invariance under deformation of the symplectic forms, they are
isomorphic~\cite[3.7]{Cieliebak handles}. We define
$FH_*(V,\, J_V,\, \phi_V)$ as $FH_*(V_c)$ for $c$ large enough. 

Let now $(J_t,\, \phi_t)$, $t\in [0,1]$ be a finite type deformation on a Stein
manifold $V$. Given $c_0$ such that all $c\ge c_0$ are regular values
of $\phi_t$ for all $t\in [0,1]$, the big Stein domains $V_{t,c}=\{\phi_t\le
c\}$, $t\in [0,1]$, $c\ge c_0$ are all diffeomorphic and endowed with
exact symplectic forms $\omega_{t,c}=-d(d\phi_t\circ J_t)$. The
boundaries $\partial V_{t,c}$ are of restricted contact type and we
can again apply Lemma 3.7 in~\cite{Cieliebak handles} to conclude that
the Floer homology groups $FH_*(V_{t,c})$ are naturally isomorphic. In
particular $FH_*(V,\, J_0,\, \phi_0)\simeq FH_*(V,\, J_1,\, \phi_1)$. 

%Let us now consider the following particular case. Assume
%$(J_0,\, \phi_0)$ and $(J_1,\, \phi_1)$ are subcritical Stein
%structures of finite type on a manifold $V$. 
%By the previous theorem, they can be connected by a finite type Stein
%deformation and therefore  $FH_*(V,\, J_0,\, \phi_0)\simeq FH_*(V,\,
%J_1,\, \phi_1)$.

\medskip 

\noindent {\it \small Proof of Theorem C.}
Cieliebak proved in~\cite{C} that, given a subcritical Stein manifold of
finite type $(\widehat N,\, J, \, \phi)$, there exists a
Stein manifold of finite type $(V,\,J_V,\,\phi_V)$ and a diffeomorphism $F:V\times
\CC\longrightarrow \widehat N$ such that: i) $J_V\times i$ and $F^*J$ are
homotopic as almost complex structures; ii) $\phi_V+|z|^2$ and
$F^*\phi$ are subcritical Morse functions
with isotopic big Stein
domains $W_c$, respectively $W'_c$. 
The functions
$\phi_V+|z|^2$ and $F^*\phi$ are associated to handle decompositions
$H_1\cup \ldots \cup H_\ell$, $H'_1\cup \ldots \cup H'_\ell$ of $W_c$
and $W'_c$ such that $\tx{index}(H_s)=\tx{index}(H'_s)$, $1\le s \le
\ell$ and the following property holds. 
Given $f\in\tx{Diff}_0(V\times \CC)$
such that $f(W_c)=W'_c$, 
the attaching maps of $f\circ H_s$, $H'_s$, $2\le s\le
\ell$ are isotopic in
$H'_1\cup\ldots\cup H'_{s-1}$ (the condition is independent of
$f$). Two such handle
decompositions are called {\it isotopic}. 
%have the same critical points (of subcritical index) 
%and have diffeomorphic large enough sublevel sets. 

Because the handle decomposition determines the isotopy class of the
associated Morse function and the two handle decompositions above are
isotopic, we infer that the 
functions $\phi_V+|z|^2$ and $F^*\phi$ are isotopic, i.e.
there exists $f\in \tx{Diff}_0(V\times
\CC)$ with $(F^*\phi)\circ f = \phi_V+|z|^2$. 
 
We can therefore apply the previous
theorem and conclude that 
the Stein structures $(J_V\times i,\, \phi_V+|z|^2)$ and
$(F^*J,\, F^*\phi)$ can be connected by a finite type Stein
deformation. It follows that
$FH_*(V\times \CC,\, F^*J,\,
F^*\phi) \simeq FH_*(V\times \CC, \, J_V\times i,\, \phi_V+|z|^2)$.
By Proposition~\ref{prop:produc with C}, the latter homology group
vanishes.
On the other hand, by invariance of Floer homology
under symplectomorphism we have $FH_*(\widehat N,\, J,\,
\phi)\simeq FH_*(V\times \CC,\, F^*J,\, F^*\phi)$. 
 \hfill{$\square$}

%This result generalizes a theorem  of Viterbo \cite{functors1}
%stating that, for any subcritical Stein manifold $\widehat N$ and any Stein
%domain $N \subset \widehat N$, the composed map $
%FH^*(\widehat N) \stackrel {c^*} 
%\longrightarrow H^*(N, \, \partial N) \stackrel {\tx{pr}}
%\longrightarrow
%H^{2n}(N, \, \partial N)$, $n =
%\frac 1 2 \dim
%\widehat N$ is
%not surjective. 

\subsection{Symplectic geometry in product manifolds}

\noindent {\it \small Proof of Theorem B.}
The statement follows readily from the existence of the 
commutative diagram given
by the second part of 
Theorem A, taking into account the isomorphism 
 $H^{2m}(M, \, \partial 
M) \otimes H^{2n}(N, \, \partial N) \stackrel \sim 
\longrightarrow H^{2m + 2n}(M 
\times N, \, \partial (M \times N))$, $2m=\dim M$,
$2n=\dim N$. The latter is  
given by the usual K\"unneth formula in singular cohomology with
coefficients in a field. 
\hfill{$\square$} 
 
\medskip 
 
\noindent {\bf Remark.} Theorem B should  be interpreted as a
stability property for the SAWC condition.

\noindent {\bf Remark.} Floer, Hofer and Viterbo~\cite{FHV} proved
the Weinstein
conjecture in a product 
$P\times \CC^\ell$, $\ell\ge 1$ with $P$ a closed symplectically
aspherical manifold. 
The Weinstein conjecture for a product $\widehat M \times \widehat N$
with $\widehat N$ subcritical Stein and $\widehat M$ the completion of
a restricted contact type manifold has been  
proved by Frauenfelder and Schlenk in~\cite{FrSc}.

\subsection{Symplectic capacities} 
The discussion below makes use of field coefficients. Let $\delta 
>0$ be small enough.
One defines (see e.g.~\cite{functors1}) the 
capacity of a compact symplectic manifold  $M$ with contact type
boundary as 
\begin{eqnarray*}
c(M) & = & \inf \{ b > 0 \, : \, FH^{m}_{]-\delta,\, b]}(M) 
\longrightarrow H^{2m}(M, \, \partial M) \tx{ is zero } \} \\
 & = & \sup \{ b > 0 \, : \, FH^{m}_{]-\delta,\, b]}(M) 
\longrightarrow H^{2m}(M, \, \partial M) \tx{ is nonzero } \} \ .
\end{eqnarray*}
Here $2m=\dim M$. The next result
is joint work with  
A.-L. Biolley, who applies it in her study of symplectic
hyperbolicity~\cite{Anne Laure}.    
\begin{prop} \label{prop:symplectic capacities} 
Let $M$, $N$ be compact symplectic manifolds with
  boundary of restricted contact type. Then   
$$c(M\times N) \le 2 \min \big( c(M), \, c(N) \big) \ . $$ 
\end{prop} 
\demo The Main Property~(\ref{inclusion fondamentale de complexes})
gives, for $\nu$ large enough and field coefficients, 
an arrow {\footnotesize$\bigoplus 
_{r+s=m+n}FH^r_{]-\delta,\frac b 2]}(H_\nu)\otimes
FH^s_{]-\delta,\frac b 2]}(K_\nu) \longleftarrow
FH^{m+n}_{]-\delta,b]}(L_\nu)$}.
%% (the direction is reversed because of
%% the cohomological setting)
Moreover, for fixed $b$ and $\nu$ large
enough the Hamiltonians $H_\nu$, $K_\nu$ and $L_\nu$ compute the
corresponding truncated cohomology groups 
of $M$, $N$ and $M\times N$. 
Like in Theorem A.b. we get the commutative diagram 
\begin{equation*}
 {\scriptsize
 \xymatrix@R=12pt{\displaystyle{
 \bigoplus _{r+s=m+n} FH^r_{]-\delta, \, \frac b 2]}(M) 
  \otimes FH^s_{]-\delta, \, \frac b 2]}(N)} \ar[d]^{c_*^{b/2}\otimes
 c_*^{b/2}} & & \ar[ll] 
  FH^{m+n}_{]-\delta, \, b]}(M \times N) \ar[d]^{c_*^b} \\  
 H^{2m}(M,\, \partial M)\otimes H^{2n}(N,\, \partial N) & &
 \ar[ll]_{\sim} H^{2m+2n}(M\times N,\, \partial (M\times N)) 
 }
 }
 \end{equation*}
Let now $b<c(M\times N)$. Then $c_*^b\neq 0$, hence 
$c_*^{b/2}\otimes c_*^{b/2}\neq 0$ and therefore $b/2 \le \min\,
\big(c(M),\, c(N)\big)$. 
\hfill{$\square$}

\end{document}